\newtheorem{thm}{Theorem}
\newtheorem{cor}{Corollary}
\newtheorem{prop}{Proposition}

\newtheorem{lem}{Lemma}

\documentclass[11pt]{article}
\usepackage{amsfonts}

\usepackage{amsmath}
\usepackage{graphicx}
\usepackage{fullpage}
\usepackage{setspace}
\newcommand{\be}{\begin{equation}}
\newcommand{\ee}{\end{equation}}
\newcommand{\ben}{\begin{equation*}}
\newcommand{\een}{\end{equation*}}
\newcommand{\ba}{\begin{aligned}}
\newcommand{\ea}{\end{aligned}}
\newcommand{\veps}{\varepsilon}
\newcommand{\rmd}{{\rm d}}
\newcommand{\rmi}{{\rm i}}

\newcommand{\dto}{\downarrow}
\newcommand{\downto}{\downarrow}
\newcommand{\upto}{\uparrow}
\newcommand{\wh}{\widehat}
\newcommand{\what}{\widehat}

\newcommand{\whA}{\widehat{A}}

\newcommand{\pibar}{\overline{\Pi}}

\newcommand{\wtilde}{\widetilde}

\newcommand{\Xbar}{X^*}
\newcommand{\Tbar}{T^*}
\newcommand{\Xstar}{\overline{X}}
\newcommand{\Tstar}{\overline{T}}

\newcommand{\Tbarb}{T^*_b}
\newcommand{\Tstarb}{\overline{T}_b}

\newcommand{\TstarYzero}{\overline{T}^Y_0}

\newcommand{\Tstarhalf}{\overline{T}_{1/2}}

\newcommand{\TstarbY}{\overline{T}_{b,Y}}
\newcommand{\TbarbY}{T^*_{b,Y}}
\newcommand{\TstarbZ}{\overline{T}_{b,Z}}

\newcommand{\R}{\mathbb{R}}

\newcommand{\ve}{\varepsilon}
\newcommand{\fb}{{1/(1-b)}}
\newcommand{\topr}{\stackrel{\mathrm{P}}{\longrightarrow}}
\newcommand{\todr}{\stackrel{\mathrm{D}}{\longrightarrow}}

\newcommand{\halmos}{$\sqcup\!\!\!\!\sqcap$}
\newcommand{\ga}{{\alpha}}
\newcommand{\gb}{{\beta}}

\numberwithin{equation}{section}
\numberwithin{thm}{section}
\numberwithin{lem}{section}
\numberwithin{prop}{section}
\numberwithin{remark}{section}
\numberwithin{cor}{section}
\numberwithin{ex}{section}

\begin{document}

\date{}

\title{Small and Large Time Stability of the Time taken for a L\'evy Process to Cross Curved  Boundaries}
\author{Philip S. Griffin and Ross A. Maller\thanks{Research partially supported by
ARC Grant DP1092502}
\\Syracuse University and Australian National University}
\maketitle

\begin{abstract}
This paper is concerned with the small time behaviour of a L\'{e}vy process $X$.
In particular, we investigate the {\it stabilities} of the
times, $\Tstarb(r)$ and  $\Tbarb(r)$, at which
$X$, started with $X_0=0$, first leaves the space-time regions
$\{(t,y)\in\R^2: y\le rt^b, t\ge 0\}$ (one-sided exit), or $\{(t,y)\in\R^2: |y|\le rt^b, t\ge 0\}$ (two-sided exit),  $0\le b<1$,
as $r\dto 0$.  Thus essentially we determine whether or not these passage times behave  like deterministic functions in the sense of different modes of convergence; specifically convergence in probability, almost surely and in $L^p$.
In many instances these are seen to be equivalent to relative stability of the process $X$ itself.  The analogous large time problem is also discussed.

\end{abstract}

\bigskip






\setcounter{equation}{0} \section{Introduction}\label{s1}

There is a strand of research, going back to \cite{BG}, and continuing most recently in \cite{BDM}, in which the local behaviour of a L\'evy process $X_t$ is compared with that of  power law functions, $t^b$, $b\ge 0$.
Here we address this question, but take a  different line, by asking for properties of the first exit time
of the process out of space-time regions bounded, either on one side or both sides, by power law functions.
Our aim is to give a very general study of the small time stability, as
the boundary level $r\to 0$, of the one-sided exit time
\be\label{tdef}
\Tstarb(r)= \inf \{t \geq 0: X_t>rt^b\},\ r\ge 0,
\ee
and the 2-sided exit time,
\be\label{tbdef}
 \Tbarb(r)= \inf \{t \geq 0: |X_t|>rt^b\},\ r\ge 0,
 \ee
when $0\le b<1$. (We adopt the
convention that the inf of the empty set is $+\infty$.)
While not the primary motivation for this paper, in Section \ref{s6} we also include results on stability for large times as the boundary level $r\to\infty$.
When $b=0$ such results form part of classical renewal theory for L\'evy processes.

The restriction of $b$ to the interval $[0,1)$ in  \eqref{tdef} and \eqref{tbdef} involves no loss of generality, since as we show below in Proposition \ref{prop3}, neither passage time can be relatively stable when $b\ge 1$.
Thus unless otherwise mentioned we keep $0\le b<1$ in what follows.
Our study will draw out similarities as well as differences between the behaviours of $ \Tstarb(r)$ and
$\Tbarb(r)$ with respect to differing modes of stability.
By {\em relative stability at 0} of $\Tstarb(r)$, we will mean
that $\Tstarb(r)/C(r)$ converges in probability to a finite nonzero constant (which by
rescaling we can take as 1), as $r\to 0$  for a finite function
$C(r)>0$. We will show that this is precisely equivalent to the {\em
positive relative stability at 0} of the process $X$, i.e., to the property that
\begin{equation}\label{relstab}
\frac{X_t}{B(t)} \topr +1,\ {\rm as}\ t\to 0,
\end{equation}
for some norming function $B(t)>0$.
The corresponding result for the two-sided exit is that $\Tbarb(r)$ is relatively stable at 0
iff $X_t$ is relatively stable at 0 in the two-sided sense, i.e., if
\begin{equation}\label{relstab2}
\frac{|X_t|}{B(t)} \topr 1,\ {\rm as}\ t\to 0,
\end{equation}
for some function $B(t)>0$.
The statements of these results are similar, and this is exploited in one direction of the proof, but the proofs in the opposite direction are completely different.

We also consider relative stability in the a.s. sense and in $L^p$. In the former case the results for the one-sided and two-sided exit times are again similar, see Theorem \ref{thm3}, and we are again able to exploit this in one direction.
In the case of stability in $L^p$, the behaviour of the two exit times is significantly different, see Theorem \ref{thm5}.  Section \ref{s3} contains a complete discussion of these results.

Given the equivalences between the relative stability of $\Tstarb(r)$ and $\Tbarb(r)$,
and the relative stability of the original process $X$, we begin Section \ref{s2} by reprising,
and where necessary extending, the properties of a relatively stable $X$.  Our main results,
related to the stability of  $\Tstarb(r)$ and $\Tbarb(r)$, are then given in Section \ref{s3}.
Proofs of these results can be found in Section \ref{s5}, together with some preliminary results which may be of independent interest. Finally Section \ref{s6} contains results in the large time setting.
We strive for, and mostly achieve, definitive (necessary and sufficient) conditions.

We conclude this section by introducing some of the notation that will be needed in the remainder of the paper.
The setting is as follows. Suppose that $X=\{X_{t}: t \geq 0 \}$,
$X_0=0$,  is a L\'{e}vy process defined on $(\Omega, {\cal F}, P)$,
with triplet $(\gamma, \sigma^2, \Pi)$, $\Pi$ being the
L\'{e}vy measure of $X$, $\gamma\in \mathbb{R}$  and $\sigma\ge
0$. Thus the characteristic function of $X$ is given by the
L\'{e}vy-Khintchine representation, $E(e^{i\theta X_{t}}) = e^{t
\Psi(\theta)}$, where
\begin{equation}\label{lrep}
\Psi(\theta) =
 \rmi\theta \gamma - \sigma^2\theta^2/2+
\int_\R (e^{\rmi\theta x}-1-
\rmi\theta x \mathbf{1}_{\{|x|\le 1\}})\Pi(\rmd x),
\ {\rm for}\  \theta \in \mathbb{R}.
\end{equation}
If $X$ is of bounded variation, then $\sigma=0$ and the
L\'{e}vy-Khintchine exponent may be expressed in the form
\begin{equation}\label{lrepbv}
\Psi(\theta) =
 \rmi\theta \rmd +
\int_\R (e^{\rmi\theta x}-1)\Pi(\rmd x),
\ {\rm for}\  \theta \in \mathbb{R},
\end{equation}
where $\rmd= \gamma -\int x \mathbf{1}_{\{|x|\le 1\}})\Pi(\rmd x)$ is called the drift of $X$.  We will sometimes include a subscript, as in for example $\rmd_X$, to make clear the process we are referring to.   $X$ is a compound Poisson process if $\sigma_X=0$, $\Pi_X(\R)<\infty$ and $\rmd_X=0$.


Let $\pibar$ and $\pibar^\pm$ denote the tails of $\Pi$, thus
\begin{equation*}
\pibar^+(x)=\Pi\{(x,\infty)\},\ \pibar^-(x)=\Pi\{(-\infty, -
x)\},\ {\rm and}\
 \pibar(x)= \pibar^+(x)+\pibar^-(x),
\end{equation*}
for $x>0$, and define kinds of Winsorised and truncated means $A(x)$ and $\nu(x)$ by
 \begin{eqnarray}\label{Adef}
A(x)&:=& \gamma +\pibar^+(1)- \pibar^-(1)
+\int_1^x \left(\pibar^+(y)- \pibar^-(y)\right)\rmd y
\nonumber\\
&=&
\gamma + x\left(\pibar^+(x)- \pibar^-(x)\right)
+\int_{1<|y|\le x} y \Pi(\rmd y)
\nonumber\\
&=:&
\nu(x)+ x\left(\pibar^+(x)- \pibar^-(x)\right),\ x>0,
 \end{eqnarray}
 where $\int_{1<|y|\le x} =-\int_{x<|y|\le 1}$ if $x<1$.
Similarly, for variances, we set
 \begin{eqnarray}\label{Udef}
U(x)&:=& \sigma^2+2\int_0^xy\pibar(y)\rmd y
\nonumber\\
&=&
\sigma^2+ x^2\pibar(x)
+\int_{0<|y|\le x} y^2 \Pi(\rmd y)
\nonumber\\
&=:&
V(x)+ x^2\pibar(x),\ x>0.
 \end{eqnarray}
Note that, because $\int_{\{|x|\le 1\}}x^2\Pi(\rmd x)<\infty$, we have
\be\label{A0}
\lim_{x\to 0} xA(x)= \lim_{x\to 0} x\nu(x)=0.
\ee




\setcounter{equation}{0}
\section{Small Time Relative Stability of $X$}\label{s2}\
Recall that $X$ is {\em relatively stable} (in probability, as
$t\to 0$), denoted $X\in RS$ at 0,  if there is a nonstochastic function $B(t)>0$ such
that
\begin{equation}\label{pmrel}
\frac{X_t}{B(t)} \topr +1,\ {\rm or}\ \frac{X_t}{B(t)} \topr -1,\ {\rm as}\ t\to 0.
\end{equation}
(We abbreviate this to $X_t/B(t)\topr \pm 1$.)
If \eqref{pmrel} holds with a $``+"$ sign we say that $X_t$ is  {\em positively relatively stable} as $t\to 0$,
denoted $X\in PRS$; a minus sign gives negative relative stability, $NRS$.
Various properties of relative stability at 0 are developed in \cite{DM2002a}.
We need only assume $B(t)>0$ for $t>0$:
$B(t)$ is not assumed {\em a priori} to be nondecreasing, but can always be taken as such.
Further properties of relative stability in probability at 0, and of the norming function $B(t)$, are
summarized in the next proposition.

\begin{prop}\label{prop1}
There is a non-stochastic function $B(t)>0$ such that
\begin{equation}\label{3.4}
{\frac{X_t}{B(t)}}\topr \pm 1,\ \mathrm{as}\ t\to  0,
\end{equation}
if and only if
\begin{equation}\label{3.5}\ {
\sigma^2=0 \quad {\rm and}\quad  \frac{A(x)}{x\pibar(x)}}\to \pm\infty,\ \mathrm{as}\ x\to  0.
\end{equation}
(The $+$ or $-$ signs should be taken together in \eqref{3.4} and
\eqref{3.5}.) If these hold, then $|A(x)|$ is slowly varying as
$x\to  0$, and $B(t)$ is regularly varying of index 1 and $B(t)\sim t|A(B(t))|$ as $t\to  0$.
Further, $B(t)$ may be chosen to be continuous and such that $t^{-b}B(t)$ is
strictly increasing  for all $0\le b<1$.

In addition, we have
\begin{equation}\label{3.4a}
{\frac{|X_t|}{B(t)}}\topr 1,\ \mathrm{as}\ t\to  0,
\end{equation}
for a non-stochastic
function $B(t)>0$, if and only if
\begin{equation}\label{3.5a}
\sigma^2=0 \quad {\rm and}\quad {\frac{|A(x)|}{x\pibar(x)}}\to \infty,\ \mathrm{as}\ x\to  0,
\end{equation}
and this is equivalent to \eqref{3.4} and \eqref{3.5} (with either the $+$ or $-$ sign). Thus
\eqref{3.4a} implies that $\lim_{t\to 0}P(X_t>0)=1$  or $\lim_{t\to 0}P(X_t<0)=1$, just as
\eqref{3.5a} implies that $A(x)$ is of constant sign for all small
$x$, that is, $A(x)>0$ for all small $x>0$ or $A(x)<0$ for all small $x>0$.

Further, the following  conditions are also each equivalent to \eqref{3.4a}:


\noindent there exist constants
$0<c_1<c_2<\infty$ and a non-stochastic
function $\wtilde B(t)>0$ such that
\begin{equation}\label{3.4aaa}
\lim_{t\to  0} P\left(c_1<\frac{|X_t|}{\wtilde B(t)}< c_2\right) \to 1;
\end{equation}

\noindent there is a  nonstochastic
function $\what B(t)>0$ such that every sequence $t_k\to
 0$ contains a subsequence $t_{k'}\to  0$ with
\begin{equation}\label{rsc}
\frac{X_{t_{k'}}}{\what B(t_{k'})}\topr c',
\end{equation}
where $c'$ is a constant with $0<|c'|<\infty$ which may depend on
the choice of subsequence. 
\end{prop}

\bigskip
\noindent {\bf Note:} \
If $\Pi(\R)=0$ then $A(x)=\gamma$ for all $x>0$, and the meaning of the limit in  \eqref{3.5} is  that $\gamma>0$ when the limit is $\infty$ and $\gamma<0$ when the limit is $-\infty$.  This corresponds to the case that $X_t=\gamma t+\sigma W_t$ is Brownian motion with drift, and it's clear that $X\in PRS$ ($X\in NRS$) at 0  iff
$\sigma^2=0$ and $\gamma>0$ ($\gamma<0$).  In this case $B(t)=|\gamma|t$.  Similarly the meaning of the limit in  \eqref{3.5a} when $\Pi(\R)=0$ is that $\gamma\neq0$.

\bigskip \noindent {\bf Proof of  Proposition \ref{prop1}}.\
See Theorem 2.2 of \cite{DM2002a} for the equivalence of \eqref{3.4} and \eqref{3.5},
and the properties of $B(\cdot)$ and $A(\cdot)$.
(A blanket assumption of $\Pi(\R)>0$ is made in \cite{DM2002a}, but it is unnecessary in any of the instances where references are made to \cite{DM2002a} in this paper.  One way to see this is to add an independent rate 1 Poisson process to $X$ and use that the resulting process agrees with $X$ at sufficiently small times.)
The strict monotonicity of
$t^{-b} B(t)$ for $0\le b<1$ follows easily from the regular variation of $B$; see for example, Section 1.5.2 of \cite{BGT}.

Clearly \eqref{3.5} implies \eqref{3.5a} and the converse holds by continuity of $A$.
Further, it is trivial that \eqref{3.4} implies \eqref{3.4a} and \eqref{3.4a} implies \eqref{3.4aaa}.
Also \eqref{3.4aaa} implies \eqref{rsc} because, under \eqref{3.4aaa}, every sequence $t_k\to
 0$ contains a subsequence $t_{k'}\to  0$ such that $X_{t_{k'}}/\wtilde{B}(t_{k'})\to Z'$ ,
where  $Z'$ is an infinitely divisible  random variable with $P(c_1\le |Z'|\le c_2)=1$.
Thus, $Z'$ is bounded a.s., hence is a constant, $c'$, say, with $|c'|\in [c_1,c_2]$. Hence we may take $\wh B = \wtilde B$ in \eqref{rsc}.
Thus to complete the proof of Proposition \ref{prop1}, it suffices to show
\eqref{rsc} implies \eqref{3.5a}.

Assume  \eqref{rsc} holds. Then every sequence $t_k\to 0$ contains a subsequence $t_{k'}\to  0$
with
\begin{equation}\label{subs}
\frac{X_{t_{k'}}}{\wh B(t_{k'})}\topr c',
\end{equation}
for some  $c'\neq 0$.  We first show that this condition holds if $\wh B$ is replaced by any function $D$   with $D(t)\in {\cal{L}}_t$ for all $t>0$, where
\ben
{\cal{L}}_t=\{{\text{limit points of $\wh B$ at $t$}}\}\cup \{\wh B(t)\}.
\een
Since $P(X_t= 0)>0$ for some $t>0$ precisely when $X$ is compound Poisson, it follows from \eqref{subs}
that $P(X_t\ne 0)=1$ for all $t>0$. Thus if $0\in {\cal{L}}_t$, then along some sequence $s\to t$,
we have $|X_s|/\wh B(s)\topr\infty$ .
From this it follows  that $0\notin {\cal{L}}_t$ if $t$ is sufficiently small.
Now take any sequence $t_k\to 0$.  Choose $s_k$ so that
\ben
 \frac{\wh B(s_k)}{D(t_k)}\to 1\quad \text{and}\quad \frac{X_{|t_k-s_k|}}{D(t_k)}\topr 0.
\een
The former is possible since $D(t_k)\in {\cal{L}}_{t_k}$, and the latter since $X_t\topr 0$ as $t\to 0$.  Now choose a subsequence $s_{k'}$ of $s_{k}$ so that
${X_{s_{k'}}}/{\wh B(s_{k'})}\topr c'$ where $c'\neq 0$.  Then
\ben
\frac{X_{t_{k'}}}{D(t_{k'})}=\frac{X_{s_{k'}}}{\wh B(s_{k'})}\ \frac{\wh B(s_{k'})}{D(t_{k'})}+ \frac{X_{t_{k'}}-X_{s_{k'}}}{D(t_{k'})}\topr c'.
\een
Thus
\be\label{Dcon}
\text{every sequence $t_k\to
 0$ contains a subsequence $t_{k'}\to  0$ with
${X_{t_{k'}}}/{D(t_{k'})}\topr c'\neq 0$.}
\ee
From the convergence criteria for infinitely divisible laws,
e.g. Theorem 15.14 of Kallenberg \cite{Kall},  this is equivalent to every sequence $t_k\to
 0$ containing a subsequence $t_{k'}\to  0$ such that for every $\ve >0$,
\begin{equation}\label{kall1}
\lim_{t_{k'}\to  0} t_{k'} \pibar(\ve D(t_{k'})) =0, \quad
\lim_{t_{k'}\to  0} \frac{t_{k'} V(D(t_{k'}))}{D^2(t_{k'})}=0,
\quad {\rm and}\quad
\lim_{t_{k'}\to  0} \frac{t_{k'} A(D(t_{k'}))}{D(t_{k'})}=c'\neq 0.
\end{equation}
From this we may conclude that,
\begin{equation}\label{kall}
\lim_{t\to  0} \frac{D(t) |A(D(t))|}{ V(D(t))}=\infty,
\quad {\rm and}\quad
\lim_{t\to  0} \frac{|A(D(t))|}{D(t)\pibar(D(t))}=\infty.
\end{equation}
Observe that $D(t) |A(D(t))|\to 0$ as $t\to 0$ by \eqref{A0}, since necessarily $D(t)\to 0$.  Hence, $\sigma^2\le V(D(t))\to 0$, proving the first condition in  \eqref{3.5a}.  Next, let
\ben
D_1(t)=\liminf_{s\to t} \wh B(s),\quad D_2(t)= \limsup_{s\to t} \wh B(s) \vee \wh B(t).
\een
Since $D_1$ and $D_2$ satisfy \eqref{Dcon}, it follows easily that
\be\label{DD}
\limsup_{t\to 0}\frac{D_2(t)}{D_1(t)}<\infty.
\ee
Now given $x>0$, let
\ben
t_x=\inf\{s:\wh B(s)\ge x\}.
\een
Then $t_x<\infty$ for sufficiently small $x$, $D_1(t_x)\le x \le D_2(t_x)$, and $t_x\to 0$ as $x\to 0$. Further
\ben\ba
|A(D_1(t_x))-A(x)|&\le |x\pibar(x)-D_1(t_x)\pibar(D_1(t_x))|+\int_{D_1(t_x)<|y|\le x} |y|\Pi(\rmd y)\\
&\le 2D_2(t_x)\pibar(D_1(t_x)),
\ea\een
hence by \eqref{kall} and \eqref{DD},
\be\label{AD}
\frac{A(x)}{A(D_1(t_x))}\to 1\quad\text{as } x\to 0.
\ee
Thus by \eqref{kall}, \eqref{DD} and \eqref{AD},
\ben
\frac{|A(x)|}{x\pibar(x)}\ge \frac{|A(D_1(t_x))|}{D_1(t_x)\pibar (D_1(t_x))}\ \frac{|A(x)|}{|A(D_1(t_x))|}\ \frac{D_1(t_x)}{D_2(t_x)}\to \infty,
\een
which proves the second condition in \eqref{3.5a}. \hfill\halmos\bigskip

\noindent {\bf Remarks:} \ (i)\
An equivalence for $PRS$, i.e., \eqref{relstab}, is \eqref{3.5} holding with a $``+"$ sign.
Then $A(x)>0$ for all small $x$. Symmetrically, for $NRS$.

(ii)\
For any family of events $A_t$,
we say that $A_t$ occur with probability approaching 1 (WPA1) as $t\to L$ if
$\lim_{t\to L}P(A_t)=1$. ($L$ may be 0 or $\infty$.)
We sometimes describe a situation like \eqref{3.4aaa}, i.e, for a stochastic function $Y_t$, $t\ge 0$, there exist constants $0<c_1< c_2<\infty$ such that
\begin{equation}\label{3.4aab}
\lim_{t\to L} P\left(c_1< Y_t< c_2\right) =1,
\end{equation}
by writing
$Y_t \asymp 1$ WPA1 as $t\to  L$.
The strict inequalities may be replaced by non-strict ones.

(iii)\
It is possible to have $A(x)\to 0$ as $x\to 0$, and also $X\in RS$ as $t\to 0$.
For example, take $\sigma^2=0$, and define
\[
\pibar^+(x)= \frac{1}{x(\log x)^2}, \ 0<x< e^{-1}, \quad \pibar^+(x)=0,\  x\ge e^{-1},
\]
and $\pibar^-(x)\equiv 0$. Then $A(x)= \gamma -1 -1/\log x$, for $x\le e^{-1}$.
Thus if $\gamma=1$ then $A(x)\to 0$ as $x\to 0$.
Further $A(x)/x\pibar(x)\to \infty$ as $x\to 0$, so that  $X\in PRS$ as $t\to 0$.
In this case, $X(t)/B(t)\topr 1$ where $B(t)=t/|\log t|$.
\bigskip

In studying $\Tstarb$ and $\Tbarb$ we will need the following corresponding maximal processes;
 \[
 \Xstar_t:= \sup_{0\le s\le t}X_s, \quad {\rm and}\quad
 \Xbar_t:= \sup_{0\le s\le t}|X_s|.
 \]

\begin{lem}\label{lem1}
Let $t_k$ be any sequence with $t_k\to 0$ as $k\to\infty$. Assume
$X_{t_k}/B(t_k) \topr a\in \R$, where $B(t_k)>0$, when $k\to \infty$. Then
\begin{equation}\label{rsc2}
{(i)}\ \frac{\Xbar_{t_k}}{B(t_k)} \topr |a|,\ {\rm and\ }(ii)\
\frac{\Xstar_{t_k}}{B(t_k)} \topr \max(a,0), \ {\rm as}\ k\to \infty.
\end{equation}
Conversely, (i) with $a\in\R$ implies $|X_{t_k}|/B(t_k) \topr |a|$ as $k\to\infty$.  Finally, as a partial converse to (ii), if $a>0$ and  ${\Xstar_{t}}/{B(t)} \topr a$ as $t\to 0$, then
${X_{t}}/{B(t)} \topr a$ as $t\to 0$.
\end{lem}

\bigskip \noindent {\bf Proof of  Lemma \ref{lem1}}.\
Assume $X_{t_k}/B(t_k)\topr a$ as $t_k\to 0$.
Use the decomposition in \cite{DM2002a}, Lemma 6.1, to write
\begin{equation}\label{[9]}
X_t =t\nu(h)+X_t^{(S,h)} +X_t^{(B,h)},\ \  t>0, h>0,
\end{equation}
where $X_t^{(B,h)}$ is the component of $X$ containing the jumps larger than $h$ in modulus, and $X^{(S,h)}$ is then  defined through \eqref{[9]}. We will use \eqref{[9]} with $t=t_k$ and $h=B(t_k)$.
As in \eqref{kall1}, $t_k\pibar(B(t_k))\to 0$ as $t_k\to 0$, so
\be
P(\sup_{0\le s\le t_k}|X_{s}^{(B,h)}|=0)=e^{-t_k\pibar(B(t_k))}\to 1.
\ee
Next,  $ X_t^{(S,h)}$ is a  mean 0 martingale with variance $tV(h)$, so by applying Doob's inequality to the submartingale $(X_t^{(S,h)})^2$, for any $\veps>0$,
\[
P\left(\sup_{0\le s\le t_k}| X_s^{(S,h)}|>\veps B(t_k)\right)
\le \frac{t_kV(B(t_k))}{\veps^2 B^2(t_k)} \to 0,
\]
using \eqref{kall1} again. A third use of \eqref{kall1} gives
$t_k\nu(B(t_k))\sim aB(t_k)$, so from \eqref{[9]},
\[
\frac{X_{t_k}^*}{B(t_k)}  =\sup_{0\le s\le t_k} \frac{|X_s|}{B(t_k)}
= \sup_{0\le s\le t_k} \frac{s|\nu(B(t_k)|}{B(t_k)} +o_p(1)  \topr |a|.
\]
Thus ({\it i}) is proved and  ({\it ii}) follows similarly.

Conversely, let ({\it i}) hold.  Then
$X_{t_k}^*/B(t_k)$ is stochastically bounded
and the inequality $P(|X_{t_k}|>xB(t_k))\le P(X_{t_k}^*>xB(t_k))$ for $x>0$ shows that
$X_{t_k}/B(t_k)$ is also stochastically bounded.
Thus every subsequence of $\{t_k\}$ contains a further subsequence $t_{k'}\to 0$ along which
$X_{t_{k'}}/B(t_{k'})\todr Z'$, for some infinitely divisible random variable
$Z'$ with $|Z'|\le |a|$.  As a bounded infinitely divisible random variable, $Z'$
is degenerate at $z'$, say. But then $X_{t_{k'}}^*/B(t_{k'})\topr |z'|$ by the converse direction already proved.  Hence by {\it{(i)}}, $|z'|=|a|$ and since this is true for all subsequences we get $|X_{t_{k}}|/B(t_{k})\topr |a|$.

Finally assume ${\Xstar_{t}}/{B(t)} \topr a$ as $t\to 0$ where $a>0$.  We may assume $B(t)$ is nondecreasing; see for example Lemma \ref{nondec}.  We first show that
\be\label{B2}
\limsup_{t\to 0}\frac{B(t)}{B(2t)}<1.
\ee
If not, there exists a sequence $t_k\to 0$ so that ${B(t_k)}/{B(2t_k)}\to 1$ as $k\to \infty$.  Thus
\be\label{2pr}
\frac{\Xstar_{t_k}}{B(t_k)}\topr a \quad\text{and}\quad \frac{\Xstar_{2t_k}}{B(t_k)}\topr a,\ \ \text{as}\ k\to\infty.
\ee
Let $\tau_k=\inf\{t:X_t>aB(t_k)/2\}$ and set
$Y_t=X_{\tau_k +t}- X_{\tau_k}$.  Then $P(\tau_k\le t_k)\ge P(\Xstar_{t_k}>aB(t_k)/2)\to 1$, and on $\{\tau_k\le t_k\}$,
\ben
\overline{Y}_{t_k}\le \Xstar_{2t_k}-X_{\tau_k}\le (\Xstar_{2t_k}-\Xstar_{t_k}) + (\Xstar_{t_k}- aB(t_k)/2).
\een
By \eqref{2pr},
\ben
\frac{\Xstar_{2t_k}-\Xstar_{t_k}}{B(t_k)} + \frac{\Xstar_{t_k}- aB(t_k)/2}{B(t_k)}\topr \frac a2,
\een
and hence
\ben
P(\Xstar_{t_k}> 3aB(t_k)/4)= P(\overline{Y}_{t_k}> 3aB(t_k)/4)\to 0,
\een
which is a contradiction.  Thus \eqref{B2} holds.  Now write
\be\label{xbd}
\Xstar_{2t}= \Xstar_t \vee (X_t+\Xstar'_t),
\ee
where $\Xstar'_t=\sup_{t\le s\le 2t}(X_s-X_t)$ is an independent copy of $\Xstar_t$.  Given any sequence $t_k\to 0$ we may choose a further subsequence $t_{k'}\to 0$ so that
\ben
\frac{B(t_{k'})}{B(2t_{k'})}\to \lambda'
\een
where necessarily $\lambda'\in [0,1)$ by \eqref{B2}. Setting $t=t_{k'}$ in \eqref{xbd}, dividing throughout by $B(2t_{k'})$  and taking limits, we see that
\be\label{sltp}
\frac{X_{t_{k'}}}{B(2t_{k'})}\topr a(1-\lambda')> 0.
\ee
Thus with $\wh B(t)= B(2t)$ in \eqref{rsc}, it follows that $X\in RS$ and since the subsequential limits in \eqref{sltp} are positive, $X\in PRS$.  Thus for some function $D(t)>0$, $X_t/D(t)\topr 1$.  Then from part ({\it ii}), $\Xstar_t/D(t)\topr 1$.  Hence $D(t)\sim aB(t)$ and the proof is complete.
\hfill\halmos\bigskip

\noindent {\bf Remark:}\
We are unsure whether a subsequential version of the converse to ({\it ii}), with $a>0$,  holds. Since it will not be needed in this paper we do not pursue it further.
\bigskip


One final result which will prove useful below is the following;

\begin{prop}\label{prop2}
Suppose $X_{t}/B(t) \topr 1$, where $B(t)>0$, when $t\to 0$, and let $Y_t(\lambda):= X_{\lambda t}/B(t)$ for $\lambda\ge 0$. Then
for every $\delta>0$, $0<\eta\le T<\infty$  and $0\le b<1$,
\ben
P(\sup_{\eta\le \lambda\le T}|\lambda^{-b}Y_t(\lambda)-\lambda^{1-b}|>\delta)\to 0\quad \text{as } t\to 0.
\een
\end{prop}

\bigskip \noindent {\bf Proof of Proposition \ref{prop2}}.\
By a result of Skorohod, see for example Theorem 15.17 of \cite{Kall}, for every $\delta >0$
\be\label{Sk}
P(\sup_{0\le \lambda\le T}|Y_t(\lambda)-\lambda|>\delta)\to 0\quad \text{as } t\to 0.
\ee
Thus
\ben
P(\sup_{\eta\le \lambda\le T}|\lambda^{-b}Y_t(\lambda)-\lambda^{1-b}|>\delta)\le P(\sup_{\eta\le \lambda\le T}|Y_t(\lambda)-\lambda|>\delta\eta^{b})\to 0
\een
by \eqref{Sk}.
\hfill\halmos




\setcounter{equation}{0}
\section{Relative Stability of $\Tstarb(r)$ and $\Tbarb(r)$ for Small Times}\label{s3}
Recall we always assume, unless explicitly stated otherwise, that $0\le b<1$.
The first two theorems concern the relative stability in probability or almost surely of
$\Tstarb(r)$ and $\Tbarb(r)$, as $r\to 0$.
These are shown to be equivalent to  positive relative stability at 0 of $X$ and  relative stability at 0 of $X$,
in the relevant mode of convergence, respectively.
Proposition \ref{prop3} demonstrates that there is no loss of generality in assuming $0\le b<1$, since $\Tstarb(r)$ and $\Tbarb(r)$ cannot be relatively stable,
in probability (and hence also a.s.), as $r\to 0$, when $b\ge 1$.

\begin{thm}\label{thm1a}
(a)\
Assume $X_t/B(t)\topr 1$ as $t\to 0$, where $B(t)>0$.  Then $B(t)/t^b$ may be chosen to be continuous and strictly increasing, in which case
$\Tstarb(r)/C(r)\topr 1$ as $r\to 0$, where $C(r)$ is the inverse to $B(t)/t^b$.

Conversely, assume
$\Tstarb(r)/C(r)\topr 1$ as $r\to 0$, where $C(r)>0$.  Then C(r) may be taken to be continuous and
strictly increasing with inverse $C^{-1}$, in which case $X_t/B(t)\topr 1$ as $t\to 0$, where $B(t)=t^bC^{-1}(t)$.

\noindent
(b) \
The same result holds if $X$ and $\Tstarb(r)$ are replaced by $|X|$ and $\Tbarb(r)$ respectively in (a).
\noindent

In either case, (a) or (b), the function
$C(r)$ is regularly varying with index $\fb$ as $r\to 0$.
\end{thm}

In the example given prior to Lemma \ref{lem1}, $X(t)/B(t)\topr 1$ where $B(t)=t/|\log t|$.  Thus $\Tstarb(r)/C(r)\topr 1$ and $\Tbarb(r)/C(r)\topr 1$  where $C(r)=\left((1-b)^{-1}r|\log r|\right)^{1/(1-b)}$.

\bigskip\noindent{\bf Remark:}\
Implicit in Theorem \ref{thm1a} we understand, is that  $\Tstarb(r)$ and $\Tbarb(r)$
are finite WPA1 as $r\to 0$, as a consequence of their relative stability when $X\in PRS$ or $X\in RS$.

\begin{cor}\label{corTT}
Assume $X_t/B(t)\topr 1$ as $t\to 0$, where $B(t)>0$. Then  $P(\Tstarb(r)=\Tbarb(r))\to 1$ as $r\to 0$.
\end{cor}

\begin{prop}\label{prop3}
Suppose $b\ge 1$. Then neither $\Tstarb(r)$ nor $\Tbarb(r)$ can be relatively stable, in probability,
as $r\to 0$.
\end{prop}

\begin{thm}\label{thm3}
(a)\ $\Tbarb(r)$ is almost surely (a.s.) relatively stable, i.e.,
$\Tbarb(r)/C(r)\to 1$, a.s., as $r\to 0$, for a finite function $C(r)>0$,
 iff $X$  has bounded variation with drift $\rm d_X \ne 0$.

\noindent
(b)\ $\Tstarb(r)$ is almost surely relatively stable, i.e.,
$\Tstarb(r)/C(r)\to 1$, a.s., as $r\to 0$, for a finite function $C(r)>0$, iff $X$  has bounded variation with drift $\rm d_X > 0$.

\noindent
In either case, (a) or (b), the function $C(r)$ may be chosen as $C(r)=(r/|\rmd_X|)^\fb$.
\end{thm}

The next theorem deals with the position of the process after exiting,
 in the setting of Theorem \ref{thm1a}.

\begin{thm}\label{thm2a}
\noindent
(a)\
Suppose $X_t/B(t) \topr 1$ as $t\to 0$, where $B(t)>0$ satisfies the regularity conditions of Theorem \ref{thm1a}. Let
$C$ be the inverse of $B(t)/t^b$.
Then, as $r\to 0$,
\begin{equation}\label{alot}
\frac{X_{\Tstarb(r)}}{r(\Tstarb(r))^b} \topr 1,\ \
\frac{X_{\Tstarb(r)}}{B(\Tstarb(r))} \topr 1,\ \
\frac{X_{\Tstarb(r)}}{B(C(r))} \topr 1,\ {\rm and}\ \
\frac{X_{\Tstarb(r)}}{r(C(r))^b} \topr 1.
\end{equation}

\noindent
(b)\
Suppose $|X_t|/B(t) \topr 1$ as $t\to 0$.
Then \eqref{alot} holds with $|X|$ and $\Tbarb(r)$ in place of $X$ and  $\Tstarb(r)$ respectively.
\end{thm}

\noindent{\bf Remark:}\ By Theorem 4.2 of \cite{DM2002a}, $X_t/B(t)\to 1$ a.s. as $t\to 0$  for some $B(t)>0$, is equivalent to $X$ having bounded variation with drift $\rmd_X>0$, and in that case $B(t)\sim \rmd_X t$.  Using this, it is then easy to see that the analogous result to Theorem \ref{thm2a} holds when $\topr$ is replaced throughout by a.s. convergence.\medskip

In the results so far, $\Tstarb(r)$ and $\Tbarb(r)$ have behaved very similarly.
This is not the case when it comes to stability in $L^p$ as
our final result 
shows.
When considering $E\Tstarb(r)$, the immediate problem arises as to whether or not the expectation is finite.  As shown in Theorem 1 of \cite{DM2004}, finiteness of
$E\Tstarb(r)$ for some (all) $r>0$ is equivalent to $X_t\to\infty$ a.s. as $t\to\infty$.  Since our aim is to study the local behaviour of $X$ for small times, imposing a large time condition  is most unnatural.  Thus, we remove the issue of finiteness of $E\Tstarb(r)$, by studying instead, $E(\Tstarb(r)\wedge \ve)$ as $r\to 0$ for small $\ve$.  Similarly for $E(\Tbarb(r)\wedge \ve)$.

\begin{thm}\label{thm5}
(a)\ Assume $X$  has bounded variation with drift $\rm d_X > 0$,
then
\begin{equation}\label{estar}
\lim_{\ve\to 0}\lim_{r\to 0} \frac{E(\Tstarb(r)\wedge \ve)}{C(r)}=1,
\end{equation}
where $C(r)=(r/\rmd_X)^\fb $.

\noindent
(b)\
Fix a function $C(r)>0$; then
 $\Tbarb(r)/C(r)\topr 1$ iff $(\Tbarb(r)\wedge \ve)/C(r)\to 1$ in $L^p$  for some (all) $p>0$ and some (all) $\ve>0$.
In particular, if $\Tbarb(r)/C(r)\topr 1$, then for every $p>0$ and $\ve>0$
\begin{equation}\label{ebar}
\lim_{r\to 0} \frac{E(\Tbarb(r)\wedge \ve)^p}{C(r)^p}=1.
\end{equation}
\end{thm}\bigskip

By standard uniform integrability arguments, see for example Theorem 4.5.2 of \cite{Dur}, \eqref{estar} implies that  ${(\Tstarb(r)\wedge \ve)}/(r/\rmd_X)^{1/(1-b)}\to 1$  in $L^p$  as $r\to 0$ then $\ve\to 0$, if  $p\le 1$.
However, unlike \eqref{ebar}, \eqref{estar}  does not extend to convergence of the $p$-th moment for $p>1$.   Nor does \eqref{estar} hold  without taking the additional limit as $\ve\to 0$.
To illustrate this, let
$
X_t=at-N_t
$
where $N_t$ is a rate one Poisson process and $a>0$.  Then $X$ has bounded variation with $\rmd_X=a$.
Clearly $\Tstarb(r)=(r/a)^{1/(1-b)}$ if $N_{(r/a)^{1/(1-b)}}=0$, while $\Tstarb(r)\ge a^{-1}$ if $N_{(r/a)^{1/(1-b)}}\ge 1$.
Thus for any $p>0$, if $\ve<a^{-1}$  and  $r$ is sufficiently small that $(r/a)^{1/(1-b)}<\ve$, then
\ben
E(\Tstarb(r)\wedge \ve)^p = (r/a)^{p/(1-b)}e^{-(r/a)^{1/(1-b)}} + \ve^p(1-e^{-(r/a)^{1/(1-b)}}).
\een
Hence, with $p=1$, we obtain
\ben
\lim_{r\to 0} \frac{E(\Tstarb(r)\wedge \ve)}{(r/\rmd_X)^{1/(1-b)} }=1+\ve,
\een
showing that the limit on $\ve\to 0$ is needed in \eqref{estar}.
If $p>1$, then
\ben
\lim_{r\to 0} \frac{E(\Tstarb(r)\wedge \ve)^p}{(r/\rmd_X)^{p/(1-b)} }=\infty,
\een
for every $\ve > 0$, so the first moment convergence in \eqref{estar} does not extend to $p$-th moment convergence for any $p>1$.

\bigskip\noindent {\bf Remarks:}\
(i)\ Although parts (a) and (b) of  Theorem \ref{thm1a}  are similar in content, the proof for  $\Tbarb$ is quite different to that for  $\Tstarb$.
To prove (a) we need to establish  {\it a priori} certain regularity properites of $C(r)$, whereas the proof of (b) relies heavily on the fact that bounded infinitely divisible distributions must be degenerate.

\noindent
(ii)\  In Theorem \ref{thm3} we deduce results for  $\Tstarb$ from those for $\Tbarb$, whereas
in Theorem \ref{thm2a}, we do the opposite.

\noindent
(iii)\  We restricted ourselves to the boundary functions $t\mapsto t^b$ in this paper for clarity of exposition, though it's clear that many of our arguments will go through for more general regularly varying or even dominated varying functions.




\setcounter{equation}{0} \section{Proofs}\label{s5}

We set out some preliminary results.
Throughout, take $0\le b<1$.  A key to proving Theorem \ref{thm1a} for $\Tstarb(r)$ is to obtain the {\it a priori} regularity of $C(r)$ contained in the following Proposition;

\begin{prop}\label{mainprop}
Suppose $\Tstarb(r)/C(r)\topr 1$ as $r\to 0$
for a finite function $C(r)>0$.  Then $C(r)$ may be chosen to be continuous and strictly increasing.
\end{prop}

We emphasize that no assumptions are being made on $C$ beyond positivity.  This creates several difficulties which could be avoided if we were to assume, for example, that $C$ is regularly varying.  Such an assumption, however, would clearly be unsatisfactory, and, as we show, unnecessary.
The main purpose of Proposition \ref{mainprop}, which is somewhat hidden in the proof of Theorem \ref{thm1a},  is that from $\Tstarb(r)/C(r)\topr 1$ we can conclude that $C(r)-C(r-)=o(C(r-))$ as $r\to 0$.  This latter condition is actually all that is needed, but proving the stronger continuity simplifies matters at several points.

The proposition
will be proved by a series of lemmas. Recalling \eqref{3.4aab} we begin with the following elementary result which we will apply below to the processes $ \Xstar, \Xbar, \Tstarb$ and $\Tbarb$:

\begin{lem}\label{nondec}
Let $W_t$ be any nonnegative, nondecreasing stochastic process with  $W_t\to 0 $ a.s. as $t\to 0$.  If
$W_t/D(t)\asymp 1$ WPA1 as $t\to 0$  for some non-stochastic function $D(t)>0$, then $D(t)\to 0$ and may be chosen to be nondecreasing.
If
$W_t/D(t)\topr 1$  as $t\to 0$  for some  function $D$, then again $D(t)\to 0$  and may be chosen to be nondecreasing.
\end{lem}

\bigskip \noindent {\bf Proof of Lemma \ref{nondec}}.\
Suppose that $P(c_1<W_t/D(t)<c_2)\to 1$ for some $0<c_1<c_2<\infty$, as $t\to 0$.  This trivially implies $D(t)\to 0$ as $t\to 0$.  To avoid pathological cases where $D(t)\to 0$ as $t\to 1$ for example, choose $t_0$ small enough that
$P(c_1<W_t/D(t)<c_2)\ge 1/2$ for all $0<t\le t_0$.  Then $0<\inf_{t\le s\le t_0}D(s)\le \sup_{t\le s\le t_0}D(s)<\infty$ for all $0<t\le t_0$.
Let $D^*(t)=\inf_{t\le s\le t_0}D(s)$ for $0<t\le t_0$. It then suffices to show
\ben
1\le \liminf_{t\to 0}\frac{D(t)}{D^*(t)}\le \limsup_{t\to 0}\frac{D(t)}{D^*(t)}\le \frac{c_2}{c_1}.
\een
Only the final inequality requires proof.  If this did not hold there would be a sequence $t_k\to 0$ with
$D(t_k)/D^*(t_k)\to a>c_2/c_1$.
Thus for some sequence $s_k\to 0$, $s_k\ge t_k$, we have $D(t_k)/D(s_k)\to a$.
But this leads to a contradiction since
\begin{equation*}
\frac{W_{t_k}}{D(t_k)}
\le \frac{W_{s_k}}{D(s_k)} \frac{D(s_k)}{D(t_k)},
\end{equation*}
and the LHS is $\ge c_1$ WPA1, whereas the RHS is
$\le c_2/a<c_1$ WPA1.
\hfill\halmos\medskip

Lemma \ref{nondec} clearly applies to  $ \Xstar$ and  $\Xbar$.
For application of Lemma \ref{nondec} to $\Tstarb$ and $\Tbarb$, note that in general
$\Tstarb(r)$ (respectively $\Tbarb(r)$) need not converge to $0$ a.s. as $r\to 0$.
However, when $\Tstarb(r)/C(r)\asymp 1$ (respectively $\Tbarb(r)/C(r)\asymp 1$) WPA1 as $r\to 0$,
almost sure convergence of $\Tstarb(r)$ and $\Tbarb(r)$ to 0 does occur.
This is because $\Tstarb(r)\downarrow \Tstar_0(0)$ a.s. as $r\downarrow 0$, and if $P(\Tstar_0(0)=0)<1$, then,
combined with $\Tstarb(r)/C(r)\asymp 1$, we would have $P(\Tstar_0(0)>c)=1$ for some $c>0$.
Hence $X_t\le 0$ for all $t$ and so $\Tstarb(r)=\infty$ for all $r>0$;
but this contradicts $\Tstarb(r)/C(r)\asymp 1$ WPA1. Similarly for $\Tbarb$.
For later reference, we note that the same argument holds if $\Tstarb$ is replaced by $\Tstar_f$, where
\begin{equation}\label{Tfdef}
 \Tstar_f(r)= \inf \{t \geq 0: X_t>rf(t)\},\ r\ge 0,
 \end{equation}
and $f$ is any function for which $f(t)>0$ for $t>0$ and $f(t)\to 0$ as $t\to 0$.  Similarly for $\Tbar_f$.

\begin{lem}\label{Tineq}
Suppose there is a (nondecreasing, without loss of generality) function $C(r)>0$ such that
\newline
(a) $\Tbarb(r)/C(r)\topr 1$ or (b) $\Tstarb(r)/C(r)\topr 1$ as $r\to 0$.
\newline
Then for every $\beta>1$, in either case,
\begin{equation}\label{csup}
\limsup_{r\to 0}\frac{C(\beta r)}{C(r)} <\infty.
\end{equation}
\end{lem}

\bigskip \noindent {\bf Proof of Lemma \ref{Tineq}}.\
(a)\
Let
\begin{equation}\label{Ysdef}
Y_s:= X_{\Tbarb(r)+s}-X_{\Tbarb(r)}, \ s\ge 0,
\end{equation}
and $\TbarbY(r)$  be the corresponding two-sided passage time, viz,
\be\label{TbY}
\TbarbY(r):= \inf\{s\ge 0: |Y_s|> rs^b\},\ r\ge 0.
\ee
Fix $\veps\in(0,1)$ so that $\xi:= 2^{1-b} ((1-\veps)/(1+\veps))^b>1$ and set
\ben\ba
A^+_r:&=\{X_{\Tbarb(r)}>0, \ Y_{\TbarbY(r)}>0, \ \frac{\Tbarb(r)}{C(r)}\in (1-\veps,1+\veps),\\
&\qquad\qquad\qquad\qquad \frac{\TbarbY(r)}{C(r)}\in (1-\veps,1+\veps),\ \Tbarb(\xi r)\ge (1-\veps)C(\xi r)\},\\
A^-_r:&=\{X_{\Tbarb(r)}<0, \ Y_{\TbarbY(r)}<0, \ \frac{\Tbarb(r)}{C(r)}\in (1-\veps,1+\veps),\\
&\qquad\qquad\qquad\qquad  \frac{\TbarbY(r)}{C(r)}\in (1-\veps,1+\veps),\ \Tbarb(\xi r)\ge (1-\veps)C(\xi r)\}.
\ea\een
Then on $A^+_r$ we have
\begin{eqnarray*}
X_{\Tbarb(r)+\TbarbY(r)}
&=&
X_{\Tbarb(r)}+Y_{\TbarbY(r)}
\nonumber\\
&\ge&
 r(\Tbarb(r))^b+ r(\TbarbY(r))^b
\nonumber\\
&>&
2r((1-\veps)C(r))^b
\nonumber\\
&=&
\xi r\left(2(1+\veps)C(r)\right)^b
\nonumber\\
&\ge&
\xi r\left(\Tbarb(r)+\TbarbY(r)\right)^b.
\end {eqnarray*}
Hence, still on $A^+_r$, we have
\be\label{xi}
(1-\veps)C(\xi r)\le \Tbarb(\xi r)\le \Tbarb(r)+\TbarbY(r) \le 2(1+\veps)C(r).
\ee
Replacing $X$ by $-X$ (which does not change $\Tbarb$ or $\TbarbY$) in this argument shows that \eqref{xi} also holds on $A^-_r$.

Since $P(A^+_r\cup A^-_r)>0$ for small $r$ (in fact, $\liminf_{r\to 0}P(A^+_r\cup A^-_r)\ge 1/2$), we have for small $r$
\[
\frac{C(\xi r)}{C(r)} \le \frac{2(1+\veps)}{1-\veps},
\]
proving that
\[
\limsup_{r\to 0}\frac{C(\xi r)}{C(r)} \le 2.
\]
Thus \eqref{csup} holds for $\beta=\xi$, and the general result holds,
for the two-sided case, by iteration and monotonicity.

(b)\
Exactly the same argument works for the one-sided case if $\Tbarb$ is replaced by $\Tstarb$ throughout, including in the definition of $Y$ in \eqref{Ysdef}, and $\TbarbY(r)$ is replaced by the corresponding one-sided exit time in \eqref{TbY}.  In fact the one-sided case is slightly simpler in that there is no need to consider the events $A_r^-$ since  $P(A^+_r)\to 1$ as $r\to 0$. \hfill\halmos

\begin{lem}\label{corin}
Suppose there is a function $C(r)>0$ such that $\Tstarb(r)/C(r)\topr 1$ as $r\to 0$.
Then, with $\lambda:=2^{1-b}$, we have
\begin{equation}\label{cin}
\liminf_{r\to 0}\frac{C(\lambda r)}{C(r)} \ge 2,
\end{equation}
and, with $\beta=\lambda^n$, $n=1,2,\ldots$, we have
\begin{equation}\label{cin2}
\liminf_{r\to 0}\frac{C(\beta r)}{C(r)} \ge \beta^\fb,
\end{equation}
or, equivalently, with $\alpha=\lambda^{-n}$, $n=1,2,\ldots$,
\begin{equation}\label{cin3}
\liminf_{r\to 0}\frac{C(r)}{C(\alpha r)} \ge \frac{1}{\alpha^\fb}.
\end{equation}
\end{lem}

\bigskip \noindent {\bf Proof of Lemma \ref{corin}}.\
Fix $\veps\in(0,1)$ and let
\begin{equation*}
Z_s:= X_{(1-\veps)C(r)+s}- X_{(1-\veps)C(r)},\ s\ge 0,
\end{equation*}
and
$\TstarbZ:=\inf\{t\ge 0: Z_t>rt^b\}$,\ $r\ge 0$.
Then on
\[
A_r:=\left\{(1-\veps)C(r)< \Tstarb(r), (1-\veps)C(r)< \TstarbZ(r), \Tstarb(\lambda r)\le (1+\veps)C(\lambda r)\right\}
\]
we claim
\begin{equation}\label{pgb0}
X_t\le 2^{1-b} rt^b,\ {\rm for\ all}\ 0\le t\le 2(1-\veps)C(r).
\end{equation}
This is trivial for $0\le t\le (1-\veps)C(r)$, while for $0< s\le (1-\veps)C(r)$, on $A_r$,
\begin{eqnarray*}
X_{(1-\veps)C(r)+s}
&=&
X_{(1-\veps)C(r)}+Z_s
\nonumber\\
&\le&
r\left((1-\veps)C(r)\right)^b+rs^b
\nonumber\\
&\le&
2^{1-b} r\left((1-\veps)C(r)+s\right)^b,
\end {eqnarray*}
where the last inequality follows from convexity of $x\mapsto x^b$, $0\le b<1$, which implies
$x^b+y^b\le 2^{1-b}(x+y)^b$, for $x,y>0$.
Thus we get \eqref{pgb0}.
So, on $A_r$, we have
$\Tstarb(\lambda r)\ge 2(1-\veps) C(r)$, where, recall, $\lambda=2^{1-b}$.
Since $P(A_r)>0$ for small $r$ (in fact $P(A_r)\to 1$ as $r\to 0$) this gives
\[
2(1-\veps)C(r)\le (1+\veps)C(\lambda r).
\]
Letting $r\to 0$ then $\ve\to 0$ yields
\eqref{cin}.

For \eqref{cin2}, let $\beta=\lambda^n$ and
write
\[
\frac{C(\beta r)}{C(r)}= \prod_{k=1}^n \frac{C(\lambda^k r)}{C(\lambda ^{k-1}r)},
\]
from which we get $\liminf_{r\to 0}C(\beta r)/C(r) \ge 2^n$.
But $2^n=(\lambda^\fb)^n=\beta^\fb$.

Finally, \eqref{cin3} follows immediately from  \eqref{cin2}.   \hfill\halmos\bigskip

Now we need a little analysis.
Fix $n\ge 1$, $a=a_n>0$ and $\alpha\in(0,1)$, and consider the following curves for $t\ge a$:
\ben
y_1=r_{n}t^b
\qquad {\rm and}\qquad
y_2=\alpha r_{n+1}(t-a)^b+r_{n+1}a^b,
\een
where for notational convenience we let $r_{n}=1/n$. (A picture which also includes the curve $y=r_{n+1}t^b$ is helpful).
We wish to estimate where these curves intersect.  For this it is more convenient to consider them in the new coordinate system:
\ben
s=t-a, \qquad
Y=y-r_{n+1}a^b,
\een
in which they become
\begin{equation}\label{fdef}
Y_1=r_{n}(s+a)^b-r_{n+1}{a^b}:=f_1(s)
\qquad {\rm and}\qquad
Y_2=\alpha r_{n+1}{s^b}:=f_2(s).
\end{equation}
Elementary calculus shows that $f_1'=f_2'$ iff
\[
s=\frac{(\alpha n)^\fb a}
{(n+1)^\fb - (\alpha n)^\fb}.
\]
Thus the curves cross at most twice. We will show that they cross at exactly 2 points and estimate the positions of these points.

To do this, first note that the function
\begin{equation}\label{gdef}
g(x):= \frac{(1+x)^b-1}{x^b}, \ x>0,
\end {equation}
is strictly increasing on $(0,\infty)$, with $g(x)\downto 0$ as $x\downto 0$,
and  $g(x)\upto 1$ as $x\upto \infty$.
For $\alpha\in(0,1)$ define $c(\alpha)=g^{-1}(\alpha)$,
and, for $c>0$,
\begin{equation}\label{Rdef}
R_{n}(c):=\frac{ca}{(\alpha n)^{1/b}}
\quad {\rm and}\quad
\wh R_{n}(c):=ca.
\end {equation}

Now we need:

\begin{lem}\label{corin2}
Define $f_1$, $f_2$, $R_{n}(c)$, and $\wh R_{n}(c)$ as in
\eqref{fdef} and \eqref{Rdef}. Fix $\alpha\in(0,1)$, but allow $a=a_n>0$ to vary with $n$.
 Then for large $n$,
\begin{equation}\label{cineq}
f_1(R_n(c)) {> \atop <} f_2(R_n(c))\ {\rm if}\ c {< \atop >} 1;
\end{equation}
and
\begin{equation}\label{c-1ineq}
f_1(\wh R_n(c)) {< \atop >} f_2(\wh R_n(c))\ {\rm if}\ c {< \atop >} c(\alpha).
\end{equation}
Consequently, for any $\veps\in (0,1)$, if $n$ is sufficiently large, then
\be\label{f2>f1}
f_2(s)>f_1(s)\ \ \text{for all}\ \ s\in \big(R_{n}(1+\veps), \wh R_{n}((1-\veps)c(\alpha))\big).
\ee
\end{lem}

\bigskip \noindent {\bf Proof of Lemma \ref{corin2}}.\
First, as $n\to\infty$,
\begin{eqnarray*}
n(n+1)f_1(R_n(c))
&=&
(n+1)\left(1+c/(\alpha n)^{1/b}\right)^b a^b-n a^b
\nonumber\\
&=&
(n+1)\left(1+bc/(\alpha n)^{1/b}+O(1/n^{2/b})\right) a^b-n a^b
\nonumber\\&\to &
a^b, \ {\rm since}\ b<1,
\end {eqnarray*}
while
\[
n(n+1)f_2(R_n(c))=\alpha n(ca/(\alpha n)^{1/b})^b
=c^ba^b,
\]
which proves the first statement.

For the second, we have that
\begin{eqnarray*}
\frac{f_1(\wh R_n(c))}{f_2(\wh R_n(c))}
&=&
\frac{(n+1)\left(1+c\right)^b a^b-n a^b}
{\alpha n (ca)^b}
\nonumber\\&=&
\frac{\left(1+c\right)^b -1} {\alpha c^b}
+\frac{\left(1+c\right)^b } {\alpha nc^b}.
\end {eqnarray*}
Since the second term on the RHS tends to 0, the result then follows from the definition of
$c(\alpha)$ and the monotonicity of $g$.

Finally, \eqref{f2>f1} follows immediately from \eqref{cineq} and \eqref{c-1ineq}.\hfill\halmos

\begin{lem}\label{cey}
Suppose $\Tstarb(r)/C(r)\topr 1$, as $r\to 0$, where $C(r)>0$. Then with  $r_n=1/n$, we have
\[
\limsup_{n\to\infty}\frac{C(r_{n})}{C(r_{n+1})}\le 1.
\]
\end{lem}

\bigskip \noindent {\bf Proof of Lemma \ref{cey}}.\
Fix $\veps\in(0,1)$ small enough that
\begin{equation}\label{pgb1}
\frac{b^\fb(1+\veps)^3}{(1-\veps)^2}<1.
\end{equation}
Recall $c(\alpha)=g^{-1}(\alpha)$, where $g$ is defined in \eqref{gdef}.
As $x\downto 0$, we have
\[
\frac{(1+x)^b-1}{x^b}\sim \frac{bx}{x^b} =bx^{1-b},
\]
thus $bc(\alpha)^{1-b}\sim \alpha$,  as $\alpha\downto 0$, or, equivalently,
\begin{equation*}
\lim_{\alpha\downto 0}\frac{c(\alpha)}{\alpha^\fb}=\frac{1}{b^\fb}.
\end{equation*}
Hence we can choose $\alpha>0$ of the form $\alpha=\lambda^{-k}$  small enough that
\begin{equation}\label{pgb2}
\frac{\alpha^\fb}{b^\fb} <(1+\veps)c(\alpha).
\end{equation}
For $n\ge 1$ let
\ben
A_{n}
=
\left\{\frac {\Tstarb(r_{n})}{C(r_{n})}\in (1-\veps, 1+\veps), \frac {\Tstarb(r_{n+1})}{C(r_{n+1})}\in (1-\veps, 1+\veps),\frac {\TstarbY(\ga r_{n+1})}{C(\ga r_{n+1})}\in (1-\veps, 1+\veps)\right\}
\een
where $Y_s= X_{\Tstarb(r_{n+1})+s}-X_{\Tstarb(r_{n+1})}, \ s\ge 0$, and $\TstarbY(r):= \inf\{s\ge 0: Y_s> rs^b\},\ r\ge 0.$
Then on $A_{n}$ we have, for $n$ sufficiently large, depending on $\alpha$,
\begin{eqnarray*}
\frac{(1+\veps)\Tstarb(r_{n+1})}{(\alpha (n+1))^{1/b}}
&\le&
\frac{(1+\veps)^2C(r_{n+1})}{(\alpha (n+1))^{1/b}}
\nonumber\\
&\le&
(1-\veps)C(\alpha r_{n+1}) \quad {\rm (by\ Lemma\ \ref{Tineq},\ and\ taking}\ n\ {\rm large\ enough)}
\nonumber\\
&<&
 \TstarbY(\alpha r_{n+1})
\nonumber\\
&<&
(1+\veps)C(\alpha r_{n+1})
\nonumber\\
&\le&
(1+\veps)^2\alpha^\fb  C(r_{n+1})\quad  {\rm (by\ \eqref{cin3})}
\nonumber\\
&\le&
\frac{(1+\veps)^2}{(1-\veps)} \alpha^\fb \Tstarb( r_{n+1})
\nonumber\\
&\le&
 \frac{(1+\veps)^3}{(1-\veps)} b^\fb c(\alpha) \Tstarb( r_{n+1})\quad {\rm (by\ \eqref{pgb2})}
\nonumber\\
&\le&
(1-\veps)c(\alpha) \Tstarb(r_{n+1}) \quad  {\rm (by\ \eqref{pgb1})}.
\end{eqnarray*}
Thus with $a=a_{n+1}=\Tstarb(r_{n+1})$ in \eqref{fdef} and \eqref{Rdef}, we have shown that on $A_{n}$, for large $n$,
\be\label{419}
R_{n+1}(1+\veps)<\TstarbY(\alpha r_{n+1})<\wh R_{n+1}((1-\veps)c(\alpha)).
\ee
This means that, on $A_{n}$,
\ben\ba
X_{\Tstarb(r_{n+1})+\TstarbY(\alpha r_{n+1})}&= X_{\Tstarb(r_{n+1})}+Y_{\TstarbY(\alpha r_{n+1})}\\
&\ge r_{n+1}{(\Tstarb(r_{n+1}))^b} + {f_2(\TstarbY(\alpha r_{n+1}))}\\
&> r_{n+1}{(\Tstarb(r_{n+1}))^b} + {f_1(\TstarbY(\alpha r_{n+1}))}
\quad  {\rm (by\ \eqref{f2>f1}\ and\ \eqref{419})}
\\
&=r_{n}{(\Tstarb(r_{n+1})+\TstarbY(\alpha r_{n+1}))^b}
\ea\een
by \eqref{fdef},
and so
\[
\Tstarb(r_{n})\le \Tstarb(r_{n+1})+\TstarbY(\alpha r_{n+1}).
\]
Since $P(A_{n})>0$ for large $n$ this implies
\[
(1-\veps)C(r_{n})\le (1+\veps)C(r_{n+1})+(1+\veps)C(\alpha r_{n+1}).
\]
Hence
\[
\limsup_{n\to\infty}\frac{C(r_{n})}{C(r_{n+1})}
\le \left(\frac{1+\veps} {1-\veps}\right) \left(1+\limsup_{n\to\infty}\frac{C(\alpha r_{n+1})}{C(r_{n+1})}\right)
\le \left(\frac{1+\veps} {1-\veps}\right) \left(1+\alpha^\fb\right),
\]
by \eqref{cin3}. Now let $\alpha\downto 0$ then $\veps\downto 0$ to complete the proof.\hfill\halmos

\bigskip \noindent {\bf Proof of Proposition \ref{mainprop}}.\
By Lemma \ref{nondec} and the paragraph following, we may assume $C$ is nondecreasing and $C(r)\to 0$ as $r\to 0$.  With $r_{n}=1/n$ as above, define
\[
D(r_{n})= C(r_{n})(1+r_n),\ \  n\ge 1,
\]
and interpolate $D(r_{n})$ linearly for $0<r\le 1$.
Then $D$ is continuous and strictly increasing on $(0,1]$. Thus to complete the proof of Proposition \ref{mainprop}, it suffices to show that
\[
\frac{D(r)}{C(r)}\to 1\ {\rm as}\ r\to 0.
\]
This follows easily from Lemma \ref{cey}, because, given $r\in (0,1]$, by letting $n$ satisfy $r_{n+1}< r\le r_{n}$,  we obtain
\[
D(r)\le D(r_{n}) = C(r_{n})(1+r_n)\le \left[C(r) +(C(r_{n})-C(r_{n+1}))\right](1+r_n),
\]
hence
\[
\frac{D(r)}{C(r)} \le \left\{1+\frac{C(r_{n})-C(r_{n+1})}{C(r_{n+1})}\right\}(1+r_n) \to 1\ {\rm as}\ n\to \infty,
\]
while
\[
D(r)\ge  D(r_{n+1}) \ge C(r_{n+1})\ge C(r) -(C(r_{n})-C(r_{n+1})),
\]
and so, also,
\[
\frac{D(r)}{C(r)} \ge 1-\frac{C(r_{n})-C(r_{n+1})}{C(r_{n+1})} \to 1\ {\rm as}\ n\to \infty.
\] \hfill\halmos

\bigskip \noindent {\bf Proof of Theorem \ref{thm1a}}.\
We first show that if $X$ is positively relatively stable, then so are $\Tstarb$ and $\Tbarb$.
This will prove one direction of ({\it a}), and also one direction of ({\it b}) because if $|X_t|/B(t)\topr  1$,
then by Proposition \ref{prop1}, either $X\in PRS$ or $X\in NRS$, and in the latter case we simply apply the result to $-X$ rather than $X$, which does not change $\Tbarb$.
Thus assume  there is a non-stochastic
function $B(t)>0$ such that
$X_t/B(t)\topr  1$, as $t\to 0$, assumed to have the properties listed in Proposition \ref{prop1}. Let $C(r)$ be
the inverse function to $B(t)/t^b$, uniquely defined because $B(t)/t^b$ is chosen continuous and strictly increasing.
Thus we have $B(C(r))=rC(r)^b$.
Since $\Tbarb(r)\le \Tstarb(r)$, it suffices to show
\be\label{1149}
P\left(\Tstarb(r)>(1+\veps)C(r)\right)\to 0\ \ {\rm and}\ \ P\left(\Tbarb(r)<(1-\veps)C(r)\right)\to 0\ \ {\rm as}\ r\to 0
\ee
for every $\ve>0$.  Now for any $\eta\in (0,1+\ve)$,
\begin{eqnarray*}
P\left(\Tstarb(r)>(1+\veps)C(r)\right)
&\le&
P\left(\sup_{0< t\le (1+\veps)C(r)} \frac{X_t}{t^b}\le r\right)
\nonumber\\
&\le&
P\left(\sup_{\eta\le \lambda\le 1+\veps}
 \frac{\lambda^{-b}X_{\lambda C(r)}}{B(C(r))}
 \le 1\right)\to 0
\end{eqnarray*}
by Proposition \ref{prop2}.
A similar argument shows that for any $\ve\in(0,1)$ and $\eta\in (0,1-\ve)$
\begin{eqnarray*}
P\left(\Tbarb(r)<(1-\veps)C(r)\right)
&\le&
P\left(\Tbarb<\eta C(r)\right) + P\left(\sup_{\eta\le \lambda\le 1-\veps}
 \frac{\lambda^{-b}|X_{\lambda C(r)}|}{B(C(r))}
 \ge 1\right).
\end{eqnarray*}
As above, the second term converges to $0$ as $r\to 0$ by Proposition \ref{prop2}.  For the first we use the L\'evy process version of Remark 2.1 and Proposition 2.1 in \cite{DG2}, which translated into our notation, and using that $A$ is slowly varying, gives, for some universal constant $c$,
\ben\ba
P\left(\Tbarb(r)<\eta C(r)\right) &\le  \frac{c\eta C(r) A(r(\eta C(r))^b)}{r(\eta C(r))^b}\\
& \sim \frac{c\eta^{1-b} C(r) A(rC(r)^b)}{rC(r)^b}\\
&\to c \eta^{1-b}
\ea\een
as $r\to 0$, since $B(C(r))=r C(r)^b$, $A(\cdot)$ is slowly varying at 0,
and  $tA(B(t))/B(t)\to 1$ as $t\to 0$, by Proposition \ref{prop1}.  Letting $\eta\to 0$ completes the proof of \eqref{1149}.

We now come to the converse direction. We first consider (a).  Thus
assume there exists a finite function $C(r)>0$ such that  $\Tstarb(r)/C(r)\topr 1$ as $r\to 0$.
Then $C(r)\to 0$ as $r\to 0$, and we may assume that $C(\cdot)$ is continuous and strictly increasing. Thus $B(t):=t^bC^{-1}(t)$ is uniquely defined and $t^{-b}B(t)\downarrow 0$ as $t\downarrow 0$.
We first show that, for each $\delta>0$,
\begin{equation}\label{pg1}
\lim_{t\to 0}P\left(\frac{X_t}{B(t)}>1+\delta\right)=0.
\end{equation}
To see this, take $t>0$ and $\lambda>0$, and define $r=C^{-1}(t/\lambda)$, so that $\lambda C(r)=t$.
On the event
$\{\Tstarb(r)/C(r)> \lambda\}$
we have
$X_s\le rs^b$ for all $0\le s\le \lambda C(r)=t$, and hence
\[
\Xstar_t\le t^b C^{-1}(t/\lambda) =\lambda^bB(t/\lambda)
\]
on that event.
Thus for every  $0<\lambda<1$,
\begin{equation}\label{pg2}
\liminf_{t\to 0}P\left(\frac{\Xstar_t}{B(t/\lambda)}\le \lambda^b \right)
\ge \liminf_{r\to 0}P\left(\frac{\Tstarb(r)}{C(r)}> \lambda\right)=1.
\end{equation}
Now given $\delta>0$, choose $0<\ve<1$ so that $(1-\veps)^b+\veps^b<1+\delta$.   This is possible since $(1-\veps)^b+\veps^b\downarrow 1$ as
$\veps\downarrow 0$.  Hence
\[
P\left(\frac{X_t}{B(t)}>1+\delta\right)
\le
P\left(\frac{\Xstar_t}{B(t)}>1+\delta\right)
\le
P\left(\frac{\Xstar_{(1-\veps)t}}{B(t)}>(1-\veps)^b\right)
+
P\left(\frac{\Xstar_{\veps t}}{B(t)}>\veps^b\right)
\to 0
\]
as $r\to 0, \ {\rm by\ \eqref{pg2}}$.

Next, still assuming that  $\Tstarb(r)/C(r)\topr 1$ as $r\to 0$,
we show that,  for each $\delta>0$,
\begin{equation}\label{pg3}
\lim_{t\to 0}P\left(\frac{X_t}{B(t)}>1-\delta\right)=1.
\end{equation}
To see this, take $\delta>0$ and choose $\veps\in(0,1)$ small enough  that
\begin{equation}\label{pg4}
(2\veps)^b\le \delta
\quad {\rm and}\quad
(1+2\veps)(1-\veps) >1.
\end{equation}
Given $t>0$, set $r=C^{-1}(t/(1-\ve))$; thus $t=(1-\veps)C(r)$.
Now since  $\Tstarb(r)/C(r)\topr 1$, we have
\[
1-\veps<\frac{\Tstarb(r)}{C(r)}<(1+2\veps)(1-\veps)
\]
WPA1 as $r\to 0$.
On this event
\[
X_s > rs^b\ {\rm for\ some}\  s\in[(1-\veps)C(r),(1+2\veps)(1-\veps)C(r)]
=[t,(1+2\veps)t].
\]
Thus for this $s$,
\be\label{s11}
X_s > rs^b\ge C^{-1}\left(\frac{t}{1-\veps}\right) t^b
\ge C^{-1}\left(t\right) t^b
=B(t).
\ee
Now suppose \eqref{pg3} fails. Then along a subsequence $t_k\to 0$, with probability bounded away from 0, we have
\[
X_s-X_{t_k} > B(t_k)-(1-\delta)B(t_k) = \delta B(t_k),\
\]
for some $ s\in[t_k,(1+2\veps)t_k]$, by \eqref{s11}.
Thus with probability bounded away from 0
\be\label{Xp}
\Xstar_{2\veps t_k}> \delta B(t_k)=p_k (2\veps t_k)^b,
\ee
where
\be\label{pt0}
p_k:= \frac{\delta B(t_k)}{(2\veps t_k)^b}\to 0.
\ee
But
\[
C^{-1}(t_k)=\frac{B(t_k)}{t_k^b} =\frac{(2\veps)^b p_k}{\delta}\le p_k
\]
by \eqref{pg4}, and so,
$t_k\le C(p_k)$.
Hence by \eqref{Xp} and \eqref{pt0}, with probability bounded away from 0, $\Tstarb(p_k)\le  2\veps t_k\le 2\veps C(p_k)$ where $p_k\to 0$.
This contradicts the relative stability of $\Tstarb(\cdot)$.
Thus \eqref{pg3} holds and together with \eqref{pg1}
this proves $X\in PRS$.

We now consider (b).
Thus suppose $\Tbarb(r)/C(r)\topr 1$ as $r\to 0$ for a function $C(r)>0$.  Then $C(r)\to 0$ as $r\to 0$, and we may assume, by Lemma \ref{nondec}, that $C(\cdot)$ is  nondecreasing. (Note that we are not assuming {\em a priori} that $C(\cdot)$ is
continuous and strictly increasing, in this proof.) For any $t>0$, define $ C^{-1}(t)=\inf\{r>0:C(r)\ge t\}$.   Then
\ben
C(C^{-1}(t)-)\le t\le C(C^{-1}(t)+),
\een
and so by Lemma \ref{Tineq}, for some constants $0<c_1\le c_2<\infty$
\be\label{dCt}
c_1\ C(C^{-1}(t))\le t\le c_2\  C(C^{-1}(t)),
\ee
if $t$ is sufficiently small.

Observe that for any $\lambda>0$,
\[
P\left(\Tbarb(r)\le \lambda C(r)\right)
\ge  P\left(\sup_{0< t\le \lambda C(r)} \frac{|X_t|}{t^b}>r\right)
\ge  P\left(\frac{\Xbar_{\lambda C(r)}}{(\lambda C(r))^b}>r\right).
\]
The LHS  tends to 0 as $r\to 0$ for $\lambda\in(0,1)$, so we have, for such $\lambda$,
\be\label{idl}
\lim_{r\to 0} P\left(\Xbar_{\lambda C(r)} \le \lambda ^b D(r)\right)
=1
\ee
where $D(r)= r C(r)^b.$
Now take any sequence $t_k\to 0$ and define $r_k= C^{-1}(t_k)$.  Note that by \eqref{dCt}, for large $k$
\be\label{rt}
c_1\ C(r_k)\le t_k\le c_2\ C(r_k).
\ee
Setting $\lambda=1/2$, it follows from \eqref{idl} that along a further subsequence $t_{k'}$,
\[
\frac{X_{ C(r_{k'})/2}}{D(r_{k'})} \todr Z'(1/2),
\]
where
$|Z'(1/2)|\le (1/2)^b$ a.s. Since $Z'(1/2)$ is infinitely  divisible, this means $Z'(1/2)$ is degenerate at a constant, $c'(1/2)$, say. By considering characteristic functions for example, this then implies
\be\label{X_C}
\frac{X_{\lambda C(r_{k'})}}{D(r_{k'})} \topr c'(\lambda),
\ee
for all $\lambda>0$, where the constants $c'(\lambda)$ satisfy $c'(\lambda)=\lambda c'(1)$.


We next show $c'(1)\neq0$.  If not, \eqref{X_C} gives
$
X_{\lambda C(r_{k'})}/D(r_{k'})\topr 0
$
for every $\lambda>0$, and so  by Lemma \ref{lem1}
$
X^*_{\lambda C(r_{k'})}/D(r_{k'})\topr 0
$
for every $\lambda>0$.
Fix $\ve\in(0,1)$. Then for any $\delta\in(0,1)$, $s,r>0$,
\begin{eqnarray}\label{thm3.1}
P\left(\Tbarb(r)\ge s\right)
&\ge&  P\left(\sup_{0<t\le s} \frac{|X_t|}{t^b}\le r\right)
\nonumber\\
&\ge&  P\left(\sup_{\delta s< t\le s} |X_t|\le r (\delta s)^b,
\sup_{0<t\le \delta s} \frac{|X_t|}{t^b}\le r\right)
\nonumber\\
&\ge&  P\left(\sup_{\delta s< t\le s} |X_t-X_{\delta s}|\le \ve r(\delta s)^b,
\sup_{0<t\le \delta s} \frac{|X_t|}{t^b}\le (1-\ve)r\right)
\nonumber\\
&\ge&  P\left(\Xbar_{(1-\delta)s}\le \ve r(\delta s)^b\right)
P\left(\Tbarb((1-\ve)r)> \delta s\right).
\end{eqnarray}
Now by Lemma \ref{Tineq}, for $r$ small enough
\[
1 \le \frac{C(r)}{C((1-\ve)r)}
\le c_\ve,
\]
for some $c_\ve<\infty$. Thus if we let
\[
\delta= \frac{1-\veps}{(1+\veps)c_\ve}, \quad \lambda=(1-\delta)(1+\ve)
\]
and substitute $r=r_k$, $s=(1+\veps)C(r_k)$  into \eqref{thm3.1},  we obtain for large $k$
\begin{equation*}
P\left(\Tbarb(r_k)\ge (1+\veps)C(r_k)\right)
\ge P\left(\Xbar_{\lambda C(r_k)}\le \ve[\delta(1+\ve)]^b D(r_k)\right)
P\left(\frac{\Tbarb((1-\ve)r_k)}{C((1-\ve)r_k)}
>1-\veps\right).
\end{equation*}
But along the sequence $k'$, the LHS converges to $0$ while both terms on the RHS converge to $1$.  Thus it must be the case that $c'(1)\neq 0$.

By again considering characteristic functions, \eqref{X_C} easily extends to
\[
\frac{X_{\lambda_{k'} C(r_{k'})}}{D(r_{k'})} \topr c'(\lambda),
\]
if $\lambda_{k'}\to \lambda\in[0,\infty)$.  By choosing a further subsequence $k''$ of $k'$ if necessary, it follows from \eqref{rt} that for some $\wh\lambda\in (0,\infty)$,
\ben
\frac{t_{k''}}{C(r_{k''})}\to \wh\lambda.
\een
Consequently
\[
\frac{X_{t_{k''}}}{D(r_{k''})} \topr c'(\wh\lambda),
\]
where $c'(\wh \lambda)=\wh\lambda c'(1)\neq0$.
Hence every sequence $t_k\to 0$ contains a subsequence $t_{k''}\to 0$ with
$X_{t_{k''}}/D(r_{k''})$ converging to a finite nonzero constant. Thus,
by Proposition \ref{prop1}, we have $X\in RS$.

Finally under (a) or (b), the proofs show that $C(r)$ may be taken as the inverse of the continuous and strictly increasing function $B(t)/t^b$ where $B(t)$ is regularly varying with index 1.  Hence
$C(r)$ is regularly varying with index $1/(1-b)$ and may be taken to be continuous and strictly increasing.
\hfill\halmos

\bigskip \noindent {\bf Proof of Corollary \ref{corTT}}.\
Assume $X_t/B(t)\topr 1$, where $B(t)/t^b$ is continuous and strictly increasing and $B(t)$ is regularly varying with index 1. Then by Theorem \ref{thm1a}
\ben
\frac{\Tbarb(r)}{C(r)}\topr 1\ \ \text{and}\ \ \frac{\Tstarb(r)}{C(r)}\topr 1,
\een
where $C$ is the inverse of $B(t)/t^b$.  In particular $B(C(r))=rC(r)^b$.
Fix $\ve\in (0,1)$.  On
\ben
A_r=\left\{\frac{\Tbarb(r)}{C(r)}\in(1-\ve,1+\ve), \frac{\Tstarb(r)}{C(r)}\in(1-\ve,1+\ve), \Tbarb(r)\neq \Tstarb(r)\right\},
\een
it must be the case that $X_{\Tbarb(r)}<0$ and so
\ben
X_{\Tbarb(r)+s} -X_{\Tbarb(r)}>2r((1-\ve)C(r))^b= 2(1-\ve)^bB(C(r))
\een
for some $0\le s\le 2\ve C(r)$.  Hence if $\liminf_{r\to 0}P(\Tstarb(r)=\Tbarb(r))< 1$, then $\limsup_{r\to 0}P(A_r)>0$ and so
\ben
\limsup_{r\to 0}P(\Xstar_{2\ve C(r)}>2(1-\ve)^bB(C(r)))>0.
\een
Using the regular variation of $B$, this contradicts $X_t/B(t) \topr 1$ if $\ve$ is sufficiently small, by Lemma \ref{lem1}.
\hfill\halmos\bigskip

\bigskip \noindent {\bf Proof of Proposition \ref{prop3}}.\
We will prove this in a little more generality than stated.
Assume $f:(0,\infty) \mapsto (0,\infty)$ is such that $f(x)\to 0$ as $x\to 0$, and there exists $\veps>0$ for which
\begin{equation}\label{b1.1}
f(x+y)\ge f(x)+f(y),\ \  {\rm for\ all}\ \  0< x\le \veps y
\end{equation}
if $y$ is sufficiently small.

For the one-sided exit, recall the definition of $\overline{T}_f(r)$ in \eqref{Tfdef},
and assume there is a $C(r)>0$ such that $\overline{T}_f(r)/C(r)\topr 1$ as $r\to 0$.
By Lemma \ref{nondec} and the paragraph following it, $C(r)\to 0$ and we may assume
$C(r)$ is nondecreasing.
Fix $\veps\in(0,1/2)$ so that \eqref{b1.1} holds.
Observe that if $r$ is sufficiently small,
\begin{equation}\label{fin}
rf(s)+rf((1-\veps/2)C(r))\le rf((1-\veps/2)C(r)+s)
\end{equation}
for all $0\le s\le \veps(1-\veps/2)C(r)$.
Since $\veps(1-\veps/2)\ge 3\veps/4$, \eqref{fin} holds in particular
for all $0\le s\le 3\veps C(r)/4$.
Now let
\[
Y_s:= X_{(1-\veps/2)C(r)+s}-X_{(1-\veps/2)C(r)},\ s\ge 0
\]
and
\[
\overline{T}^Y_f(r):= \inf\{s\ge 0: Y_s>rf(s)\},\ r\ge 0.
\]
Then by \eqref{fin}
\begin{equation}\label{Tfin}
P\left(\overline{T}_f(r) \le (1+\veps/4)C(r)\right)
\le
P\left(\overline{T}_f(r)\le(1-\veps/2)C(r)\right)+P\left(\overline{T}^{Y}_f(r)\le 3\veps C(r)/4\right).
\end{equation}
Since the LHS of \eqref{Tfin} tends to 1,
while $P\left(\overline{T}_f(r)\le(1-\veps/2)C(r)\right)$ tends to 0,
we may conclude that
\[
\lim_{r\to 0}P\left(\overline{T}_f(r) \le 3\veps C(r)/4\right)= \lim_{r\to 0}P\left(\overline{T}^{Y}_f(r) \le 3\veps C(r)/4\right)=1,
\]
which is a contradiction, since $3\veps/4<1$.

The proof for the 2-sided exit is virtually the same; simply replace
$\overline{T}^Y_f(r)$ by $T^{*,Y}_f(r):= \inf\{s>0: |Y_s|>rf(s)\}$, $r\ge 0$.
\eqref{Tfin} holds with this replacement.
\hfill\halmos\bigskip

If $f(x)=x^b$ with $b\ge 1$ it's easy to check that $f$ satisfies \eqref{b1.1}.
More generally, if, for small $y$, $f'$ is increasing and
\[
f'(y)\ge \frac{f(x)}{x}\ {\rm for}\ 0< x\le \veps y,
\]
then \eqref{b1.1} holds. For example, $f(x)=x/|\log x|$ satisfies this condition.
If \eqref{b1.1} holds for all $x,y$ and $\ve=1$, then $f$ is superadditive, so the
 proposition holds for this class of functions also.\medskip

Before proceeding to the proof of Theorem \ref{thm3}, we need some preliminary results which may be of independent interest.
We begin with a corollary to a result of Erickson \cite{E}.  In it we allow for the possibility of a killed subordinator $Y,$ that is, a process obtained from a proper subordinator $\cal Y$ by killing at an independent exponential time $e(q)$ with mean $q^{-1}$;  thus
\ben
Y_t=\begin{cases} {\cal Y}_t, & \text {if}\ t<e(q),\\
\partial, & \text {if}\ t\ge e(q),
\end{cases}
\een
where $\partial$ is a cemetery state.  The extension from proper to killed subordinators is trivial, but is needed below.

\begin{prop}\label{gaps} Let $Y$ be a (possibly killed) subordinator, then
\ben
\lim_{t\to 0}\frac{Y_{t-}}{Y_t}=1\  \text{a.s.}\quad\text{iff}\quad \rmd_Y>0.
\een
\end{prop}

\bigskip \noindent {\bf Proof of Proposition \ref{gaps}}.\
Since killing does not affect the drift, it  suffices to prove the result for proper subordinators. If $\rmd_Y>0$ then $Y_{t}/t\to  \rmd_Y$ a.s.
by Proposition 3.8 of \cite{Bert}.
This is easily seen to imply $Y_{t-}/t\to  \rmd_Y$ a.s., and so
\be\label{Y}
\lim_{t\to 0}\frac{Y_{t-}}{Y_t}=1\  \text{a.s.}
\ee

Conversely, assume \eqref{Y} holds  and by way of contradiction assume $\rmd_Y=0$.  Clearly \eqref{Y} implies that $Y$ can not be compound Poisson.
Since $\rmd_Y=0$, and $\sigma_Y=0$ since $Y$ is a subordinator, it must be the case that $\Pi(\R)=\infty$.  We may now apply Theorem 2 of \cite{E}.  In the terminology of \cite{E}, under \eqref{Y}, the function $h(x)=x$ is not a small gap function, and hence
\be\label{EC}
\int_0^1 \frac{x\Pi(\rmd x)}{\int_0^x\pibar(y)\rmd y}<\infty;
\ee
see Theorem 2 and the first paragraph of page 459 in \cite{E}.  Thus
\ben
\frac{\int_0^z x\Pi(\rmd x)}{\int_0^z\pibar(y)\rmd y}\le \int_0^z \frac{x\Pi(\rmd x)}{\int_0^x\pibar(y)\rmd y} \to 0\ {\rm as}\ z\to 0.
\een
Since
\ben
\int_0^x\pibar(y)\rmd y=x\pibar(x)+\int_0^x y\Pi(\rmd y),
\een
\eqref{EC} then implies that
\be
\int_0^1 \frac{\Pi(\rmd x)}{\pibar(x)}<\infty.
\ee
But this is a contradiction since $\Pi(\R)=\infty$.
\hfill\halmos

\bigskip\noindent{\bf Remark:}\
Let $\Delta Y_t=Y_t-Y_{t-}$.  Then Proposition \ref{gaps} can be rephrased as
\ben
\lim_{t\to 0}\frac{\Delta Y_t}{Y_{t-}}=0\  \text{a.s.}\quad{\it iff}\quad \rmd_Y>0.
\een
By a similar argument, one can check that if $Y$ is not a compound Poisson subordinator, then
\ben
\limsup_{t\to 0}\frac{\Delta Y_t}{Y_{t-}}=\infty\  \text{a.s.}\quad{\it iff}\quad \rmd_Y=0.
\een
\medskip

At this point we need to introduce a little fluctuation theory for which we refer  to \cite{Bert}, \cite{doneystf} or \cite{kypbook}.
Let  $L_t$ denote the local time of $X$ at its maximum and $(L^{-1}_t,H_t)_{t \geq 0}$
the bivariate ascending ladder process of $X$.  If $X_t\to -\infty$ a.s then  $(L^{-1},H)$ may be obtained from a proper bivariate subordinator by exponential killing.  Let  $\kappa(\cdot,\cdot)$ denote the Laplace exponent of $(L^{-1},H)$. Then
\begin{eqnarray} \label{kapexp}
\kappa(\ga,\gb)
= k+\rmd_{L^{-1}}\ga+\rmd_H\gb+\int_{t\ge 0}\int_{h\ge 0}
\left(1-e^{-\ga t-\gb h}\right)
\Pi_{{ L}^{-1}, { H}}(\rmd t, \rmd h),\quad \ga, \gb \ge 0,
\end{eqnarray}
where $\rmd_{L^{-1}}\ge 0$ and $\rmd_H\ge 0$ are drift constants, $\Pi_{{ L}^{-1}, { H}}$ is the L\'evy measure of $(L^{-1},H)$ and $k\ge 0$ is the killing rate.
\bigskip

\begin{lem}\label{BV} For any L\'evy process $X$,
\ben
\text{$X$ is of bounded variation and $\rmd_X>0$ iff $\rmd_{L^{-1}}>0$ and $\rmd_H>0$}.
\een
In that case $\rmd_X=\rmd_H/\rmd_{L^{-1}}$.
\end{lem}

\bigskip \noindent {\bf Proof of Lemma \ref{BV}}.\
By Theorem 2.2{\it (b)(ii)} of \cite{GMstab}, $X$ is of bounded variation and $\rmd_X>0$ iff $\sigma=0$, $\rmd_{L^{-1}}>0$ and $\rmd_H>0$, in which case  $\rmd_X=\rmd_H/\rmd_{L^{-1}}$.
Thus it suffices to show that if $\rmd _{L^{-1}}>0$ and $\rmd _{H}>0$, then $\sigma=0$.  If $X$ is compound Poisson then so is $H$, and consequently $\rmd_H=0$. But $\rmd _{H}>0$, and so $X$ can not be compound Poisson. Thus by Corollary 4.4({\it v}) of \cite{doneystf},  since $\rmd _{L^{-1}}>0$, the downward ladder height process $\wh H$ is compound Poisson.  Hence $\rmd_{\wh H}=0$, and so by Corollary 4.4({\it i}) of \cite{doneystf}, $\sigma=0$.
\hfill\halmos\bigskip

The following result is a companion to Theorem 4.2 in \cite{DM2002a}.

\begin{prop}\label{rsXbar} Assume that
\be\label{Xbar1}
\lim_{t\to 0}\frac{\Xstar_t}{B(t)}=1\  \text{a.s.}
\ee
for some function $B(t)>0$, then $X$ is of bounded variation with $\rmd_X>0$.
\end{prop}

\bigskip \noindent {\bf Proof of Proposition \ref{rsXbar}}.\
Since ${\Xstar_t}/{B(t)}\topr 1$ as $t\to 0$, it follows from Lemma \ref{lem1} that ${X_t}/{B(t)}\topr 1$ as $t\to 0$.  Hence by Proposition \ref{prop1} we may assume that
$B$ is increasing, continuous and regularly varying with index $1$. By continuity of $B$ it follows easily from \eqref{Xbar1} that
\be\label{Xbar2}
\lim_{t\to 0}\frac{\Xstar_{t-}}{B(t)}=1 \  \text{a.s.}
\ee
Fix $t<L_\infty$. If $L^{-1}_{t-}<L^{-1}_{t}$, then $\Xstar_{L^{-1}_{t-}}= \Xstar_{(L^{-1}_{t})-}$ since $\Xstar$ does not increase on the interval $(L^{-1}_{t-}, L^{-1}_{t})$.
Hence by \eqref{Xbar1} and \eqref{Xbar2},
\be\label{B1}
\frac{B(L^{-1}_{t-})}{B(L^{-1}_{t})}=I(L^{-1}_{t-}=L^{-1}_{t})+
I(L^{-1}_{t-}<L^{-1}_{t})\frac{B(L^{-1}_{t-})}{\Xstar_{L^{-1}_{t-}}}\
\frac{\Xstar_{(L^{-1}_{t})-}}{B(L^{-1}_{t})} \to 1\ \text{a.s.}
\ee
Consequently, by regular variation of $B$,
\ben
\lim_{t\to 0}\frac{L^{-1}_{t-}}{L^{-1}_{t}}=1\  \text{a.s.}
\een
Thus $\rmd _{L^{-1}}>0$ by Proposition \ref{gaps}.  Further
\ben
\lim_{t\to 0}\frac{H_{t-}}{H_{t}}=\lim_{t\to 0}\frac{\Xstar_{(L^{-1}_{t-})-}}{\Xstar_{L^{-1}_{t}}}= 1\  \text{a.s.}
\een
by \eqref{Xbar1}-\eqref{B1}.  Hence $\rmd_H>0$, again by Proposition \ref{gaps}.
The result now follows from Lemma \ref{BV}.
\hfill\halmos\bigskip

\bigskip \noindent {\bf Proof of Theorem \ref{thm3}}.\
(a)\
Assume that
$X$ has bounded variation with drift $\rmd_X\neq 0$.  
Then  by Theorem 4.2 of \cite{DM2002a}
\be\label{asd}
 \lim_{t\to 0}\frac{X_t}{t}= \rmd_X\ \ \rm{ a.s.}
\ee
This is easily seen to imply $t^{-1}X_{t-}\to  \rmd_X$ a.s., and so
\be\label{asD}
\lim_{t\to 0}\frac{\Delta X_t}t=0\ \ \rm{ a.s.}
\ee
where $\Delta X_t=X_t-X_{t-}$.
Now
\ben
r(\Tbarb(r))^b\le |X_{\Tbarb(r)}|\le r(\Tbarb(r))^b + |\Delta X_{\Tbarb(r)}|,
\een
so dividing through by  $\Tbarb(r)$ and using \eqref{asd} and \eqref{asD} we obtain
\begin{equation}\label{as2}
 \lim_{r\to 0}r(\Tbarb(r))^{b-1}=|\rmd_X|\ \rm{ a.s.,}
 \end{equation}
or equivalently
 \begin{equation}\label{as3}
 \lim_{r\to 0}\frac{\Tbarb(r)}{r^\fb}=\frac{1}{|\rmd_X|^\fb}\ {\rm a.s.}
 \end{equation}

 Conversely, assume $\Tbarb(r)/C(r)\to 1$ a.s. as $r\to 0$.  Then by Theorem \ref{thm1a}, we may assume $C$ is regularly varying with index $1/(1-b)$, continuous, strictly increasing and $C(r)\to 0$ as $r\to 0$.  Let $B(t)=t^bC^{-1}(t)$ and fix $\eta<1<\lambda$.  Then a.s. for small $t$, we have $\eta t<\Tbarb(C^{-1}(t))<\lambda t$.  For such $t$
 \be\label{UB}
 X^*_{\lambda t}\ge X^*_{\Tbarb(C^{-1}(t))}\ge C^{-1}(t)(\Tbarb(C^{-1}(t)))^b > C^{-1}(t)(\eta t)^b = \eta^b B(t)
 \ee
 and
 \be\label{LB}
 X^*_{\eta t}\le X^*_{\Tbarb(C^{-1}(t))-}\le C^{-1}(t)(\Tbarb(C^{-1}(t)))^b <  C^{-1}(t)(\lambda t)^b =  \lambda^b B(t).
 \ee
Since $B(t)$ is regularly varying with index 1 as $t\to 0$, it then easily follows from \eqref{UB} and \eqref{LB} that
$\lim_{t\to 0}\Xbar_t/B(t) =1$ a.s.  An examination of the proof of Theorem 4.2 in \cite{DM2002a}, shows that from this we may conclude that
$X$ has bounded variation with $\rmd_X\neq 0$.

(b)\
Now assume that
$X$ has bounded variation with drift $\rmd_X> 0$.
Then $\lim_{t\to0}t^{-1}X_t=\rmd_X> 0$ a.s., and so $P(\Tstarb(r)=\Tbarb(r) \text{ for all small  } r)=1$.  Consequently, from part (a),
 \begin{equation}\label{as5}
 \lim_{r\to 0}\frac{\Tstarb(r)}{r^\fb}=\frac{1}{\rmd_X^{\ \fb}}\ {\rm a.s.}
 \end{equation}

The proof of the converse for $\Tstarb(r)$ is virtually the same as for $\Tbarb(r)$.  Arguing as above, first show $\lim_{t\to 0}\Xstar_t/B(t) =1$ a.s., and then use Proposition \ref{rsXbar} to complete the proof. \hfill\halmos

\bigskip \noindent {\bf Proof of Theorem \ref{thm2a}}.\

 (a)\
 Assume the first condition fails.  Fix $\xi>1$  so that, with
 \[
 A_r= \left\{X_{\Tstarb(r)}>\xi r(\Tstarb(r))^b\right\},
\]
we have $\limsup_{r\to 0}P(A_r)>0$.  Now $\Tstarb(r)=\Tstarb(\xi r)$ on $A_r$, and by Theorem \ref{thm1a}, for arbitrary $\delta\in(0,1)$,
$\Tstarb(r)\le (1+\delta)C(r)$ and  $\Tstarb(\xi r)\ge (1-\delta)C(\xi r)$ hold WPA1 as $r\to 0$, hence
\be
\limsup_{r\to 0}P\left(A_r\cap\{\Tstarb(r)\le (1+\delta)C(r), \Tstarb(\xi r)\ge (1-\delta)C(\xi r)\}\right)>0.
\ee
Thus
\[
\liminf_{r\to 0}\frac{C(\xi r)}{C(r)} \le \frac{1+\delta}{1-\delta}.
\]
Letting $\delta\dto 0$ we get a contradiction, since by Theorem \ref{thm1a}, $C(r)$ is regularly varying with index $1/(1-b)$.
Hence
\[
 \frac{X_{\Tstarb(r)}}{r(\Tstarb(r))^b}\topr 1, \ {\rm as}\ r\to 0.
 \]
Then, since $B$ is regularly varying with index 1, $B(C(r))=rC(r)^b$ and  $\Tstarb(r)/C(r)\topr 1$,
the remaining relationships in Theorem \ref{thm2a} follow.

(b)\
Assume $|X_t|/B(t) \topr 1$ as $t\to 0$. By Proposition \ref{prop1}, there is no loss of generality in assuming
$X_t/B(t) \topr 1$.  The result then follows from Corollary \ref{corTT} and part (a).
 \hfill\halmos

\bigskip \noindent {\bf Proof of Theorem \ref{thm5}}.\
(a)\ Assume that $X$ has bounded variation with drift $\rmd_X> 0$.
Then by Theorem \ref{thm3}, for every $\ve >0$,
\[
\frac{\Tstarb(r)\wedge \ve}{(r/\rmd_X)^\fb}\to 1 \ \text{ a.s. as } r\to 0.
\]
Hence by Fatou's Lemma,
\begin{equation}\label{low}
\liminf_{r\to 0} \frac{E(\Tstarb(r)\wedge \ve)}{(r/\rmd_X)^\fb}\ge 1.
\ee
Letting $\ve\to 0$ proves a lower bound for \eqref{estar}.

For the upper bound, we first prove the result for $b=0$.  Recall the
bivariate ascending ladder process  $(L^{-1},H)$ of $X$, and its Laplace exponent $\kappa(\cdot,\cdot)$ given by
\eqref{kapexp}. Clearly
\ben
\lim_{q\to\infty}\frac{\kappa(q,0)}{q}= \rmd_{L^{-1}},
\een
and by Lemma \ref{BV}, we have $\rmd_{L^{-1}}> 0$,  $\rmd_H> 0$ and  $\rmd_X=\rmd_H/\rmd_{L^{-1}}$.
For $q>0$, let $e(q)$ be independent of $X$ and have exponential distribution with mean $q^{-1}$.  Then a straightforward calculation shows that for any $\ve\in (0,\infty]$
 (with the obvious interpretation when $\ve=\infty$),
\be\ba\label{Tveq}
E(\Tstar_0(r)\wedge \ve \wedge e(q))= \int_0^\infty e^{-qs}P(\Tstar_0(r)\wedge \ve >s)\rmd s = q^{-1} P(\Tstar_0(r)\wedge \ve > e(q)).
\ea\ee
Hence by (8) on p.174 of \cite{Bert},
\ben\ba
E(\Tstar_0(r)\wedge e(q))= \frac{\kappa(q,0)}{q}V^q(r),
\ea\een
where
\ben
V^q(r)=\int_0^\infty E(e^{-qL^{-1}_t};H_t\le r)\rmd t.
\een
 Now  the Laplace transform of $V^q$ is given by
\ben
\wh V^q(\lambda):=\lambda\int_0^\infty e^{-\lambda r}V^q(r)\rmd r =\frac 1{\kappa(q,\lambda)},\ \ \lambda>0,
\een
and so
\ben
\lambda\wh V^q(\lambda)=\frac{\lambda}{\kappa(q,\lambda)}\to\frac 1{\rmd_H}\ \ \rm{as}\ \lambda\to\infty.
\een
Thus
\ben
\frac{V^q(r)}{r}\to\frac 1{\rmd_H}\ \ {\rm{as}}\ r \to 0,
\een
by Karamata's Tauberian Theorem; see Theorem 1.7.1$'$ in \cite{BGT}.  Hence
\be\label{qlim}
\lim_{q\to\infty}\lim_{r\to 0} \frac{E(\Tstar_0(r)\wedge e(q))}{r}
=\lim_{q\to\infty}\lim_{r\to 0} \frac{\kappa(q,0)V^q(r)}{qr}
=\frac{\rmd_{L^{-1}}}{\rmd_H}=\frac 1{\rmd_X}.
\ee
  Next fix $\ve>0$ and $q\in (0,\infty)$.  Then
\ben\ba
E(\Tstar_0(r)\wedge \ve)&=E(\Tstar_0(r)\wedge \ve; \Tstar_0(r)\wedge\ve\le e(q))+E(\Tstar_0(r)\wedge \ve; \Tstar_0(r)\wedge \ve >  e(q))\\
&\le E(\Tstar_0(r)\wedge e(q)) + \ve P(\Tstar_0(r)\wedge \ve > e(q))\\
&= E(\Tstar_0(r)\wedge e(q)) + \ve q E(\Tstar_0(r)\wedge \ve \wedge e(q)) \quad ({\rm{by }}\ \eqref{Tveq})\\
&\le (1+\ve q)E(\Tstar_0(r)\wedge e(q)).
\ea\een
Thus for every $q\in (0,\infty)$,
\ben
\lim_{\ve\to 0}\lim_{r\to 0} \frac{E(\Tstar_0(r)\wedge \ve)}{r}
\le
\lim_{r\to 0} \frac{E(\Tstar_0(r)\wedge e(q))}{r}.
\een
Letting $q\to\infty$ and using \eqref{qlim}, proves the upper bound in \eqref{estar} for $b=0$.

To deal with the upper bound when $0<b<1$, introduce
\[
Y_t=\frac{X_t}{\rmd_X}-bt,\ t\ge 0.
\]
Then $Y$ has bounded variation with drift $\rmd_Y=1-b>0$. Take $r>0$ and let $\lambda_r=r/\rmd_X$.
Consider the function
\[
f(t):=\lambda_rt^b-bt, \ t\ge 0.
\]
This increases from 0 at $t=0$  to a maximum of $\lambda_r^\fb(1-b)$ at $t=\lambda_r^\fb$,
then decreases to 0 at $t=(\lambda_r/b)^\fb$.
Hence it is non-negative for $t\in [0,(\lambda_r/b)^\fb]$ and lies entirely below the horizontal line of height
$\lambda_r^\fb(1-b)$ for $t\ge 0$. Thus with
$\TstarYzero(\lambda)=\inf\{t\ge 0:Y_t>\lambda\}$,
we have
\[
\Tstarb(r)\le \TstarYzero(\lambda_r^\fb(1-b)).
\]
Thus invoking the $b=0$ result just proved, for $Y$, we have
\ben\ba
\lim_{\ve\to 0}\lim_{r\to 0} \frac{E(\Tstar_0(r)\wedge \ve)}{(r/\rmd_X)^{1/(1-b)}}
&\le (1-b)\lim_{\ve\to 0}\lim_{r\to 0} \frac{E(\TstarYzero(\lambda_r^\fb(1-b))\wedge \ve)}{\lambda_r^\fb(1-b)}\\
&=\frac{1-b}{\rmd_Y}=1.
\ea\een

(b)\
Assume $\Tbarb(r)/C(r)\topr 1$, and fix $p>0$ and $\ve >0$.  By Proposition \ref{prop1} and Theorem \ref{thm1a}, $A(r)$ is slowly varying at $0$, and we may assume without loss of generality that $A(r)>0$ for small $r$.
It then follows from \cite{p} (see also (4.3) of \cite{DG2}), when translated to the current notation, that for every $r>0,t>0$ and $m\ge 1$,
\be\label{PrU}
P(\Tbarb(r)\ge t)\le P(\Tbar_0(rt^b)\ge t)\le \left(\frac{cmr}{t^{1-b}A(rt^b)}\right)^m,
\ee
where $c\in (0,\infty)$ denotes an unimportant constant that may change from one usage to the next.
Furthermore, again by Proposition \ref{prop1} and Theorem \ref{thm1a}, we may assume that $C^{-1}(t)=t^{-b}B(t)$
where $tA(B(t))/B(t)\to 1$ as $t\to 0$.  Setting $t=C(r)$, and using $B(C(r))=r C(r)^b$, this gives
\be\label{ACr}
\lim_{r\to 0}\frac{C(r)A(rC(r)^b)}{rC(r)^b}= 1.
\ee

Choose $\xi>0$ sufficiently small that $1-b-b\xi>0$.  Since $A$ is slowly varying at $0$, there is a function $\wh A$ such that $A(r)\sim\wh A(r)$ as $r\to 0$, and $r^\xi\whA(r)$ is increasing on $(0, a]$ for some $a\in (0,\ve]$.
For any $p>0$, write
\be\ba\label{ET*1}
E(\Tbarb(r)\wedge \ve)^p&=\int_0^\ve pt^{p-1}P(\Tbarb(r)\ge t)\rmd t\\
&\le \int_0^{C(r)} pt^{p-1}\rmd t + \int_{C(r)}^a pt^{p-1}\left(\frac{cmr}{t^{1-b}A(rt^b)}\right)^m\rmd t
+ \int_a^\ve p  t^{p-1} P(\Tbarb(r)\ge a)\rmd t\\
&= I+II+III.\\
\ea\ee
Clearly $I= C(r)^p$, while by \eqref{PrU}
\ben
III\le  \ve^p\left(\frac{cmr}{a^{1-b}A(ra^b)}\right)^m
\sim \ve^p\left(\frac{cmr}{a^{1-b}A(r)}\right)^m =o(r^{(1-\eta)m}), \quad {\rm as}\ r\to 0,
\een
for any $\eta>0$, since $A$ is slowly varying.  Now $C(r)$ is regularly varying with exponent $1/(1-b)$, thus by choosing $m$ sufficiently large we see that $III=o(C(r)^p)$.  Finally for $II$ we observe that if
$C(r)\le t\le a$, then
\ben
t^{b\xi}\wh A(rt^b)\ge  C(r)^{b\xi}\wh A(rC(r)^b).
\een
Hence
\ben\ba
t^{p-1}\left(\frac{r}{t^{1-b}\wh A(rt^b)}\right)^m &\le t^{p-1-m(1-b-b\xi)}C(r)^{m(1-b-b\xi)}\left(\frac{rC(r)^b}{C(r)\wh A(rC(r)^b)}\right)^m.
\ea\een
Recalling that $1-b-b\xi>0$, we see that if $m$ is sufficiently large that $m(1-b-b\xi)>p$, then
\ben
\int_{C(r)}^at^{p-1}\left(\frac{r}{t^{1-b}\wh A(rt^b)}\right)^m \rmd t\le cC(r)^p
\een
for small $r$ by \eqref{ACr}.  Hence by \eqref{ET*1}
\ben
\limsup_{r\to 0} \frac{E(\Tbarb(r)\wedge \ve)^p}{C(r)^p}<\infty.
\een
Thus $(\Tbarb(r)\wedge \ve)/C(r)$ is bounded in $L^p$ for every $p>0$, which together with $\Tbarb(r)/C(r)\topr 1$ proves convergence in $L^p$.
The converse direction is trivial.
\hfill\halmos\bigskip

\noindent {\bf Remark:} \
If $X_t\to\infty$ a.s. as $t\to\infty$, then $E\Tstar_0(r)<\infty $ for some (every) $r>0$. In that case, if in addition $EL_1^{-1}<\infty$ and ${d_H}>0$, then by Theorem 2.3 of \cite{GMstab}
\begin{equation*}
\lim_{r\to 0} \frac{E\Tstar_0(r)}{r}=\frac{EL_1^{-1}}{d_H}.
\end{equation*}
If $0\le b<1/2$ then $E\Tbarb(r)^p<\infty $ for every $r>0$, $p>0$, and the proof of \eqref{ebar} can be modified to show
\begin{equation*}\label{ebar1}
\lim_{r\to 0} \frac{E\Tbarb(r)^p}{C(r)^p}=1.
\end{equation*}




\setcounter{equation}{0}
\section{Relative Stability of $\Tstarb(r)$ and $\Tbarb(r)$ for Large Times}\label{s6}

In this section we briefly summarise relative stability of $\Tstarb(r)$ and $\Tbarb(r)$ for large times.  All proofs are omitted.  In many cases they parallel the proofs given for small times although in some cases there are nontrivial differences.  We must first discuss the definitions of $\Tstarb(r)$ and $\Tbarb(r)$.
It is possible for $X$ to cross the $t^b$  boundary for small $t$,
but not for large $t$. For example, when $\sigma^2>0$ we have by \cite{M4} that
\[
\limsup_{t\downarrow 0} \frac{X_t}{\sqrt{2\sigma^2 t\log|\log t|}} =1\ {\rm a.s.},
\]
thus $\limsup_{t\downarrow 0}X_t/\sqrt{t}=+\infty$ a.s., while we can have in addition that $X_t$ drifts to $-\infty$ a.s. as $t\to\infty$.
If we took the infimum in \eqref{tdef} over all $t>0$, we may have that
$\Tstarhalf(r)$ is finite, in fact, takes value 0, for all $r>0$,
even though $\limsup_{t\to\infty} X_t/\sqrt{t}<\infty$ a.s.
Since we are interested in the behaviour of $X$ for large $t$, we prevent this kind of behaviour by taking the inf in \eqref{tdef}
over $t\ge 1$. Thus we define
\begin{equation}\label{5tdef}
 \Tstarb(r)= \inf \{t \geq 1: X_t>rt^b\},\ r\ge 0,
 \end{equation}
and
\begin{equation}\label{5tbdef}
 \Tbarb(r)= \inf \{t \geq 1: |X_t|>rt^b\},\ r\ge 0.
 \end{equation}

We always assume, unless explicitly stated otherwise, that $0\le b<1$.
The results up to Theorem \ref{5thm5} below, parallel those for small times if $r\to 0$ is replaced by $r\to\infty$ and the drift $\rmd_X$ is replaced by the mean $EX_1$ at the appropriate points.

\begin{thm}\label{5thm1a}
(a)\
Assume $X_t/B(t)\topr 1$ as $t\to\infty$, where $B(t)>0$.  Then $B(t)/t^b$ may be chosen to be continuous and strictly increasing, in which case
$\Tstarb(r)/C(r)\topr 1$ as $r\to\infty$, where $C(r)$ is the inverse to $B(t)/t^b$.

Conversely, assume
$\Tstarb(r)/C(r)\topr 1$ as $r\to\infty$, where $C(r)>0$.  Then C(r) may be taken to be continuous and
strictly increasing with inverse $C^{-1}$, in which case $X_t/B(t)\topr 1$ as $t\to\infty$, where $B(t)=t^bC^{-1}(t)$.

\noindent
(b) \
The same result holds if $X$  and $\Tstarb(r)$ are replaced by $|X|$ and $\Tbarb(r)$ respectively in (a).
\noindent

In either case, (a) or (b), the function
$C(r)$ is regularly varying with index $\fb$ as $r\to\infty$.
\end{thm}

\begin{cor}\label{5corTT}
Assume $X_t/B(t)\topr 1$ as $t\to\infty$, where $B(t)>0$. Then  $P(\Tstarb(r)=\Tbarb(r))\to 1$ as $r\to\infty$.
\end{cor}

\begin{prop}\label{5prop3}
Suppose $b\ge 1$. Then neither $\Tstarb(r)$ nor $\Tbarb(r)$ can be relatively stable, in probability,
as $r\to\infty$.
\end{prop}

\begin{thm}\label{5thm3}
(a)\ $\Tbarb(r)$ is almost surely (a.s.) relatively stable, i.e.,
$\Tbarb(r)/C(r)\to 1$, a.s., as $r\to\infty$, for a finite function $C(r)>0$,
 iff $E|X_1|<\infty$ and $\mu:= EX_1\ne 0$.

\noindent
(b)\ $\Tstarb(r)$ is almost surely relatively stable, i.e.,
$\Tstarb(r)/C(r)\to 1$, a.s., as $r\to\infty$, for a finite function $C(r)>0$, iff $E|X_1|<\infty$ and $\mu= EX_1>0$.

\noindent
In either case, (a) or (b), the function $C(r)$ may be chosen as $C(r)=(r/|\mu|)^\fb$.
\end{thm}

\begin{thm}\label{5thm2a}
\noindent
(a)\
Suppose $X_t/B(t) \topr 1$ as $t\to \infty$, where $B(t)>0$ satisfies the regularity conditions of Theorem \ref{5thm1a}. Let
$C$ be the inverse of $B(t)/t^b$.
Then, as $r\to \infty$,
\begin{equation}\label{5alot}
\frac{X_{\Tstarb(r)}}{r(\Tstarb(r))^b} \topr 1,\ \
\frac{X_{\Tstarb(r)}}{B(\Tstarb(r))} \topr 1,\ \
\frac{X_{\Tstarb(r)}}{B(C(r))} \topr 1,\ {\rm and}\ \
\frac{X_{\Tstarb(r)}}{r(C(r))^b} \topr 1.
\end{equation}

\noindent
(b)\
Suppose $|X_t|/B(t) \topr 1$ as $t\to \infty$.
Then \eqref{5alot} holds with $|X|$ and $\Tbarb(r)$ in place of $X$ and $\Tstarb(r)$ respectively.
\end{thm}

\noindent{\bf Remark:}\ The analogous result to Theorem \ref{5thm2a} holds when $\topr$ is replaced by a.s. convergence throughout.\medskip

For the large time version of Theorem \ref{thm5}, it is no longer appropriate to truncate the passage times since $C(r)\to \infty$ as $r\to\infty$.


\begin{thm}\label{5thm5}
(a)\ Suppose $E|X_1|<\infty$ and $\mu=EX_1>0$.
Then
\begin{equation}\label{5estar}
\lim_{r\to\infty} \frac{E\Tstarb(r)}{C(r)}=1.
\end{equation}
where $C(r)=(r/\mu)^\fb$.

\noindent
(b)\
Fix a function C(r); then
 $\Tbarb(r)/C(r)\topr 1$ iff $\Tbarb(r)/C(r)\to 1$ in $L^p$  for some (all) $p>0$.
In particular, if $\Tbarb(r)/C(r)\topr 1$, then for every $p>0$
\begin{equation}\label{5ebar}
\lim_{r\to\infty} \frac{E\Tbarb(r)^p}{C(r)^p}=1.
\end{equation}
\end{thm}
\medskip

\bigskip\noindent{\bf Remarks:}\
(i)\ Suppose $EX_1^2<\infty$ and $EX_1=0$. Then $X$ cannot be relatively stable,
so neither can
$\Tstarb(r)$ or $\Tbarb(r)$.  If $\pibar(x_0)=0$ for some $x_0>0$, then
$EX_1^2<\infty$, and so relative stability of $\Tstarb(r)$ or $\Tbarb(r)$ hinges on whether $EX_1=0$ or not.

\noindent
(ii)\  In addition to covering  one-sided passage times and also dealing with the important case $b=1/2$, omitted in  \cite{DG1}, \cite{DG2}, \cite{DM2000}, Theorem \ref{5thm1a} and associated results provide a more general approach than that of \cite{DM2000}, where the norming functions are assumed {\em a priori} to have strong regularity properties, such as regular variation, whereas we make no such assumptions.

\noindent
(iii)\  Siegmund \cite{S1} contains
the random walk version of \eqref{5estar}.
He mentions extensions of his result and possible applications to sequential confidence intervals and
hypothesis tests. We expect that similar extensions can be carried out in the general L\'evy case.


\begin{thebibliography}{99}
\footnotesize

\bibitem{Bert}
J. Bertoin,   {\em L\'{e}vy Processes}.
Cambridge University Press, Cambridge, (1996).

\bibitem{BDM}
J. Bertoin,   R.A. Doney, R.A. Maller,
Passage of L\'evy processes across power law
boundaries at small times. {\em Ann. Probab.},  36, (2008), 160--197.

\bibitem{BGT}
N.H. Bingham, C.M. Goldie,  J.L.  Teugels,
{\em Regular Variation.}
Cambridge University Press, Cambridge, (1987).

\bibitem{BG}
R.M. Blumenthal,   R.K.  Getoor,
Sample functions of stochastic processes with stationary independent increments.
{\em J. Math. Mech.}, 10, (1961), 492--516.

\bibitem{doneystf}
R.A. Doney,
{\em Fluctuation Theory for L\'evy Processes.}
Lecture Notes in Math., { 1897}, Springer, (2005).

 \bibitem{DG1}
R.A.  Doney, P.S. Griffin,   Overshoots over curved boundaries.{\em  Adv. Appl. Probab.}, 35, (2003), 417--448.

\bibitem{DG2}
R.A.  Doney, P.S. Griffin,  Overshoots over curved boundaries II. {\em Adv. Appl. Probab.}, 36, (2004), 1138--1174.

 \bibitem{DM2000}
R.A. Doney, R.A. Maller,   Random walks crossing curved boundaries: functional limit theorems, stability and asymptotic distributions for exit times and positions. {\em Adv. Appl. Probab.}, 32, (2000), 1117--1149.

\bibitem{DM2002a}
R.A. Doney, R.A. Maller,
Stability and attraction to normality for L\'{e}vy processes at zero and infinity.
{\em J. Theoretical Probab.}, { 15}, (2002), 751--792.

\bibitem{DM3}
R.A. Doney, R.A. Maller,
Stability of the overshoot for L\'evy processes.
{\em Ann. Probab.}, 30, (2002), 188--212.

\bibitem{DM2004}
R.A. Doney, R.A. Maller,
Moments of passage times for L\'evy processes.
{\em Annales de l'Institut Henri Poincar\'e (B), Probability and
Statistics}, 40, (2004), 279--297.


\bibitem{Dur}
R. Durrett,
{\em Probability: Theory and Examples}. Brooks/Cole, USA, (2005).


\bibitem{E}
K.B. Erickson,
Gaps in the range of nearly increasing processes with stationary independent increments.
{\em Z. Wahrscheinlichkeitstheorie verw. Gebiete}, 62, (1983), 449--463.






\bibitem{GMstab}
P.S. Griffin,  R.A. Maller,   Stability of the exit time for L\'evy processes.
Adv. Appl. Probab., (2011), (to appear).


\bibitem{Kall}
O. Kallenberg,
{\em Foundations of Modern Probability}. Springer, Berlin Heidelberg New York, (2001).







\bibitem{KM1999}
H. Kesten,   R.A.  Maller,
Stability and other limit laws for exit times
of random walks from a strip or a halfplane.
{\em Annales de l'Institut Henri Poincar\'e (B), Probability and
Statistics},  35, (1999), 685--734.

\bibitem{kypbook}
A. Kyprianou,
{\em Introductory Lectures on Fluctuations of L\'{e}vy
Processes with Applications}. Springer, Berlin Heidelberg New York, (2006).




\bibitem{M4}
R.A. Maller,
Small-time versions of Strassen's law for L\'evy processes.
{\em Proc. London Math. Soc.}, 98,  (2008), 531--558.


\bibitem{p} W.E. Pruitt,    The growth of random walks and L\'{e}vy
processes. {\em Ann. Probab.}, 9, (1981), 948--956.



\bibitem{S1}
D.O. Siegmund,
Some one-sided stopping rules. {\em Ann. Math. Statist.}, 38, (1967), 1641--1646.

\end{thebibliography}
\end{document}